\documentclass{article}
\usepackage{amssymb,amsmath, amsthm, amsfonts}
\usepackage{graphicx}
\usepackage{hyperref}
\usepackage{color}
\definecolor{darkgreen}{rgb}{0,0.6,0}
\definecolor{lightblue}{rgb}{.4,.68,.74}

\newcommand{\Z}{\mathbb{Z}}
\newcommand{\Q}{\mathbb{Q}}
\newcommand{\1}{{\bf 1}}
\newcommand{\D}{D\mathcal{H}}
\newcommand{\UC}{U\mathcal{C}}
\newcommand{\BNC}{B\mathcal{C}}
\newcommand{\fcd}{F^{\textrm{cd}}}

\newcommand{\twed}{\frac{1}{2}\textrm{-twisted}}

\newtheorem{theorem}{Theorem}[section]

\theoremstyle{plain}
\newtheorem{thm}[theorem]{Theorem}
\theoremstyle{plain}
\newtheorem{cor}[theorem]{Corollary}
\theoremstyle{plain}
\newtheorem{prop}[theorem]{Proposition}
\theoremstyle{plain}
\newtheorem{lemma}[theorem]{Lemma}
\theoremstyle{definition}
\newtheorem{defn}[theorem]{Definition}
\theoremstyle{definition}
\newtheorem{remark}[theorem]{Remark}
\theoremstyle{definition}
\newtheorem{example}[theorem]{Example}
\theoremstyle{definition}
\newtheorem{conj}[theorem]{Conjecture}
\theoremstyle{definition}

\begin{document}

\begin{center}
{\huge A Diagramless Link Homology}

\bigskip

Adam McDougall\\
University of Iowa
\end{center}

\bigskip

\begin{abstract}
A homology theory is defined for equivalence classes of links under isotopy in $S^3$.  Chain modules for a link $L$ are generated by certain surfaces whose boundary is $L$, using surface signature as the homological grading.  In the end, the diagramless homology of a link is found to be equal to some number of copies of the Khovanov homology of that link.  There is also a discussion of how one would generalize the diagramless homology theory (hence the theory of Khovanov homology) to links in arbitrary closed oriented 3-manifolds.
\end{abstract}

\section{Introduction}

Link homologies are typically dependent upon link diagrams.  In Khovanov homology for example (see \cite{Kh1}), one builds chain modules from a given projection of a link.  To be a link invariant, the homology needs to be invariant under the three Reidemeister moves -- something that needs proof.  In this paper, a homology theory is defined from links directly, rather than from link diagrams.  Given a link $L$, we find that the diagramless homology built from $L$ consists of some number of copies of the Khovanov homology for $L$.

Sections 2 through 4 are where the initial (complete) homology theory is built up and the necessary proofs are given.  In Section 2 preliminary definitions are given.  The goal of this section is to define what a $D^k$-surface is; $D^k$-surface is the name given to surfaces (with $k$ `crosscuts') that satisfy certain conditions.  Section 3 defines the chain modules $\mathcal{C}_{i,j,k,b}$.  Given a link $L$, the chain modules are generated by $D^k$-surfaces which have $L$ as boundary.  The homological grading $I$ is given by the signature of the surface.  Most of Section 4 is spent proving that the differential $d$ is well-defined.

Sections 5 and 6 reduce and refine the chain complexes in various ways.  Section 5 focuses on defining a Frobenius extension and using it to reduce the chain complex via `skein relations'.  The goal of Section 6 is to prove that the particular choice of ordering of the crosscuts on a $D^k$-surface is not important; different choices of an ordering for the crosscuts always yields the same homology.  The resulting (reduced) homology is called the `diagramless homology' of a link.

In Section 7, an injective chain map from the chain complex for Khovanov homology into the chain complex for the diagramless homology is given.  This process involves defining a special type of $D^k$-surface called a `state surface'.  In Section 8 we see that state surfaces actually span the entire chain complex, eventually implying that the diagramless homology of a link is equal to some number of copies of Khovanov homology.

As with any new theory, it is good to have examples.  Section 9 gives examples of how one can calculate the diagramless homology of links.  Specifically, the diagramless homology of the unknot is given for $k=0,1$ and $2$ crosscuts.

In Section 10, we discuss the possibility of generalizing this theory of links in $S^3$ to links in closed oriented 3-manifolds.  Although examples would be more difficult to compute, the existence of such a theory is evident.  A guideline is given for how to alter the definitions given in Section 2 in order to obtain the groundwork for a theory in a closed oriented 3-manifold $M$.

Some additional remarks are given in Section 11.

\section{Definitions}
\label{2}

Readers unfamiliar with the next two definitions can find more information in \cite{G-L}.

\begin{defn}
Given a surface $F$, the {\it Goeritz matrix}\footnote{The reader may be used to the definition of the Goeritz matrix being associated with generators for $H_1(F,\Z)$ as opposed to $H_1(F,\Q)$ (as in \cite{G-L}).  For this paper we work over $\Q$ for simplicity, without affecting any of our results.}  $G_F$ of $F$ is the $n \times n$ matrix whose $(i,j)$ entry is lk$(a_i,\tau a_j)$, where the $a_i$ are generators for $H_1(F,\Q)$, `lk' denotes linking number, and $\tau a_j$ is the pushoff of $2 \alpha_j$ into the complement of $F$.
\end{defn}

\begin{defn}
The {\it signature of a surface} $F$, denoted by sig$(F)$, is defined to be the signature of the Goeritz matrix $G_F$ of that surface.  Recall that the signature of a matrix is the number of positive eigenvalues minus the number of negative eigenvalues.
\end{defn}

Given a link $L$, we will only be interested in surfaces $F$ such that $\partial(F) = L$.  Our surfaces will be decorated by dots and \textit{crosscuts}, a term which is defined next.

\begin{defn}
Given a surface $F$ with boundary, a {\it crosscut} $c \in F$ is a properly embedded arc in $F$, i.e. $\partial c = c \cap \partial F$.  In other words there exists an embedding $f:[0,1] \rightarrow F$ with $f(\{0,1\}) = f([0,1]) \cap \partial F$. 
\end{defn}

\begin{remark}
In this paper, crosscuts will be given an orientation (direction).
\end{remark}

\begin{defn}
Crosscuts will be labeled as {\it active} or {\it inactive}.  \textcolor{darkgreen}{Active} crosscuts will be denoted by a \textcolor{darkgreen}{green} color, and \textcolor{red}{inactive} crosscuts by a \textcolor{red}{red} color.
\end{defn}

\begin{tabular}{ll}
	\begin{tabular}{l}
	Crosscuts on a surface are present to keep track of twisting\\ 
	that may occur in the surface.  Just as there are two different\\ 
	ways to twist a surface (left-hand \& right-hand twists), there\\ 
	are two different types of crosscuts: \textcolor{darkgreen}{active} and \textcolor{red}{inactive}.\\
	\ \\
	An example of a surface with crosscuts is given to the right $\rightarrow$
	\ \\ \ \\ \ \\ \ \\ \ \\ \ \\ \ 
	\end{tabular} 
&
	\includegraphics[scale=2]{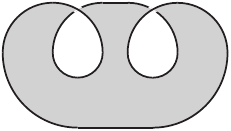}
	\thicklines
	\color{red}
	\put(-101.7,61.4){\vector(-1,4){3}}
	\color{darkgreen}
	\put(-58.5,61.3){\vector(-1,4){3}}
	\color{black}
\end{tabular} 
{\vskip -67 pt}
{\vskip -67 pt}

\begin{defn}
\label{2.6}
For a surface $F$ with crosscuts, the {\it cross-dual} of $F$, denoted $F^{\textrm{cd}}$, is the surface obtained by replacing each neighborhood of each crosscut in $F$ by the corresponding piece of (locally oriented) surface as explained below:
\begin{itemize}
\item Cut along each \textcolor{red}{inactive} crosscut and insert a locally oriented piece of surface with a left-handed $\frac{1}{2}$-twist so that the local orientation agrees with the (cut) crosscut sites.  The orientation information is kept local, and the crosscut information is forgotten.
\vspace{0pt}
	\begin{center} 
	\includegraphics[trim = 0mm 0mm 0mm 8mm, clip,scale=2]{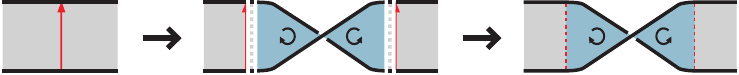}
	\put(-395,46.5){$F$}
	\put(-70,46.5){$\fcd$}
	\end{center}
\item Cut along each \textcolor{darkgreen}{active} crosscut and insert a locally oriented piece of surface with a right-handed $\frac{1}{2}$-twist so that the local orientation agrees with the (cut) crosscut sites.  The orientation information is kept local, and the crosscut information is forgotten.
	\begin{center}
	\includegraphics[trim = 0mm 8mm 0mm 0mm, clip,scale=2]{newcd.pdf}
	\put(-395,45){$F$}
	\put(-70,45){$\fcd$}
	\end{center}
\end{itemize}
\end{defn}

\begin{remark}
\label{2.7}
At times in this paper we will need to refer to the locally oriented pieces of surface inserted at (cut) crosscut sites on a cross-dual surface.  Denote the piece of (locally oriented) surface inserted along the crosscut $c$ by $F^{\textrm{cd}}|_c$.
\end{remark}

The next definition makes use of this new notation.

\begin{defn}
Given a cross-dual surface $\fcd$ with pieces of (locally oriented) surface $\{ \fcd|_{c_1}, ... ,\fcd|_{c_k} \}$, the \textit{skeleton} of $\fcd$, denoted skel($\fcd$), equals those pieces of surface along with the boundary of the cross-dual.  That is, 
$$
\textrm{skel}(\fcd) = \partial(\fcd) \cup \left( \fcd|_{c_1} \cup \cdots \cup \fcd|_{c_k} \right).
$$
\end{defn}

\begin{example}
\label{firstex}
Below is an example of a surface $F$, its cross-dual $F^{\textrm{cd}}$, and the skeleton skel($\fcd$).

{\vskip 15 pt}
\begin{center}
\includegraphics[trim = 78.4mm 0mm 21mm 0mm,clip,scale=2.5]{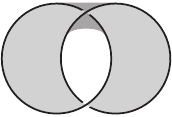}
		\thicklines
		\color{red}
		\put(-56.5,67){\vector(-1,0){18}}
		\color{darkgreen}
		\put(-38.5,3.5){\vector(-1,1){15.5}}
		\color{black}
\put(-.7,38){$\stackrel{\textrm{cross-dual}}{\longrightarrow}$}
\put(-69,86){$F$}
\ \ \ \ \ \ \ \ \ \ \ 
\includegraphics[trim = 0mm 0mm 17.5mm 0mm,clip,scale=2.5]{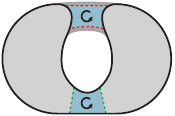}
\put(2.8,38){$\stackrel{\textrm{skeleton}}{\longrightarrow}$}
\put(-69,86){$\fcd$}
\ \ \ \ \ \ \ \ \ \ 
\includegraphics[trim = 17.5mm 0mm 0mm 0mm,clip,scale=2.5]{cd1.pdf}
\put(-80.5,86){skel($\fcd$)}
\end{center}
\end{example}

\begin{defn}
For $F$ a surface with crosscuts $\{c_1, ...,c_k\}$, refer to the components of $F - \{N(c_1), ..., N(c_k)\}$ as the {\it facets of} $F$, where $N(c)$ is a small open neighborhood in $F$ of the crosscut $c$.  On the other hand, we refer to the components of $\fcd - \textrm{skel}(\fcd)$ as the {\it facets of} $\fcd$.
\end{defn}

The following lemma tells us that the two definitions of facets above describe the same things.  This lemma will be used in later sections.

\begin{lemma}
\label{facet_lemma}
The facets of $F$ are isotopic to the corresponding facets of $F^{\textrm{cd}}$.
\end{lemma}

\emph{Proof.}  Facets of $F$ and $\fcd$ are each obtained by removing surface near crosscuts.  Since $F$ and $\fcd$ only differ near crosscuts, this means that the facets of $F$ are isotopic to the corresponding facets of $\fcd$.  The figures below show the relationship between facets of $F$ and facets of $\fcd$ near a crosscut.

\begin{center}
	\includegraphics[scale=2]{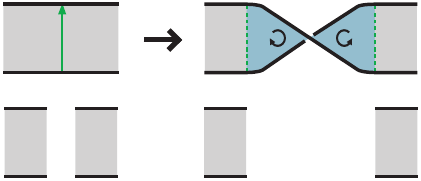}
	\put(-253,79){$F$}
	\put(2,79){$F^{\textrm{cd}}$}
	\put(-291,19){facets of $F$}
	\put(2,19){facets of $F^{\textrm{cd}}$}
	\put(-157.5,17){{\huge $\cong$}}
\end{center} \vspace{-8 pt} \qed

\begin{defn}
A {\it $D^k$-surface} is a compact surface $F$ with $k$ crosscuts $\{c_1, ...,c_k\}$ such that
\begin{itemize}
\item the crosscuts are oriented and ordered,
\item the facets of $F$ are allowed to be decorated by dots (which are \textit{not} allowed to move from one facet to another), 
\item the cross-dual $F^{\textrm{cd}}$ is orientable (this global orientatibility is independent of the local orientations of the cross-dual), and
\item there exists embedded 3-balls $B^3_+, B^3_- \subseteq S^3$ and an embedded oriented 2-sphere $\Sigma \subseteq S^3$, with $\Sigma = B^3_+ \cap B^3_-$, such that 
\begin{itemize}
\item[$\circ$] skel$(\fcd) \subseteq \Sigma$,
\item[$\circ$] all of the locally oriented pieces of surface of $\fcd$ agree with the orientation of $\Sigma$, and
\item[$\circ$] $\fcd - \textrm{skel}(\fcd) = \{\textrm{the facets of }\fcd\} \subseteq B^3_+$.
\end{itemize} 
\end{itemize}
\end{defn}

Although the crosscuts of a $D^k$-surface are to be ordered, we will define a homology theory (in Section \ref{dhsection}) for which the particular choice of crosscut ordering does not matter.

\begin{example}
It is easy to check that the surface $F$ given in Example \ref{firstex} is a $D^k$-surface (except that an ordering of the crosscuts is not given).  $\fcd$ is an orientable surface, the skeleton of $\fcd$ is planar with both pieces of surface having a positive orientation in the plane, and $\fcd$ can be viewed as having all of its facets sitting behind the plane.  Hence, $F$ is a $D^2$-surface for the Hopf link because $F$ has $k=2$ crosscuts and 
$\partial(F) = $ \ \ \ \ \ .
\put(-19,0){$\bigcirc$}
\color{white}
\put(-12.8,4.4){{\tiny $\bullet$}}
\color{black}
\put(-13,0){$\bigcirc$}
\color{white}
\put(-12.8,-1.7){{\tiny $\bullet$}}
\color{black}
\put(-14.73,-4.29){\rotatebox{47.9}{{\scriptsize \scalebox{1}[1.1]{$\smile$}}}}
\put(-18.03,-4.52){ \scalebox{1}[1.08]{\rotatebox{47.9}{{\scriptsize \scalebox{1}[1.1]{$\smile$}}}} }
\end{example}

\section{Chain Modules}
\label{3}

In the definition of $D^k$-surface, it was noted that the facets of $F$ are allowed to be decorated by dots.  Denote the total number of dots on the surface $F$ by $\delta$:
$$
\delta(F) = \# \textrm{ dots on } F.
$$

The homology has four gradings, which are each fixed by the differential $d$.  Given a $D^k$-surface $F$, define

\begin{itemize}
\item $I(F) := \textrm{sig}(F)$.  This will be our homological grading.
\item $J(F):=-\chi (F) - \textrm{sig}(F) + 2 \cdot \delta(F)$.  This will be our polynomial grading.\footnote{The $J$-grading is called the \textit{polynomial grading} due to its similarity to the polynomial grading found in Khovanov homology.}
\item $K(F) := k$.  This is the number of crosscuts on $F$.
\item $B(F) := \textrm{sig}(F) \ + $ (\# of \textcolor{darkgreen}{active} crosscuts on $F$).
\end{itemize}

The chain modules for homology are defined next.  In later sections, reduced chain modules will be introduced to produce a more interesting homology.

\begin{defn}
Given a link $L$, we let $\mathcal{C}_{i,j,k,b}(L)$ be the free module of isotopy classes of $D^k$-surfaces in $S^3$ with $I=i, J=j, K=k,B=b,$ and $\partial(F) = L$.  
\end{defn}

\section{The Differential}

Now the process of defining our differential $d:\mathcal{C}_{i,j,k,b} \rightarrow \mathcal{C}_{i+1,j,k,b}$ begins.  We define $d$ in parts, acting locally on neighborhoods of \textcolor{darkgreen}{active} crosscuts on $D^k$-surfaces in $\mathcal{C}_{i,j,k,b}$.

\begin{defn}
\label{4.1}
Given a $D^k$-surface $F$ with an \textcolor{darkgreen}{active} (and oriented) crosscut $c$ on $F$, define {\bf $d_c$} to be the map which replaces a neighborhood of $c$ in $F$ with the piece of surface shown below\footnote{For figures depicting surfaces, darker shading indicated the presence of more layers of surface.  For example, in Definition \ref{4.1}, the surface labeled $d_c(F)$ is a $\twed$ band with a second $\twed$ band attached to it.  The second $\twed$ band has a darker shading to indicate that it is in front of the other band from our viewpoint.}:
	\begin{center} 
	\includegraphics[scale=1.5]{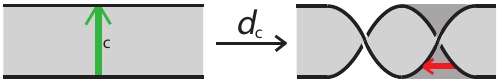}
	\put(-27.5,.5){{\footnotesize \textbf{c} }}
	\put(-228,14){$F$} 
	\put(3,14){$d_c(F)$}
	\end{center}
\end{defn}

There is a potential issue with the above definition that needs to be checked.  Given a neighborhood of a crosscut $c$ in $F$, in general there is no way of distinguishing one side of this piece of surface from the other.  One might worry that viewing a neighborhood of $c$ from the `front' side versus the `back' side will result in different surfaces after applying $d_c$.  However, this is not the case.  This is proved below.

\begin{thm}
The same surface is obtained after applying $d_c$ to either side of a neighborhood of the crosscut $c$.
\end{thm}

\emph{Proof.}
\begin{center}
\begin{tabular}{ll}
	\begin{tabular}{l}
	Make note that when the $x, y,$ and $z$ axes are referred to,\\
	the coordinate system shown to the right is being used. \\
	\ 
	\end{tabular} 
&
	$\begin{array}{l}\textit{}
		\includegraphics[trim = 0mm 0mm 27mm 0mm, clip,scale=.4]{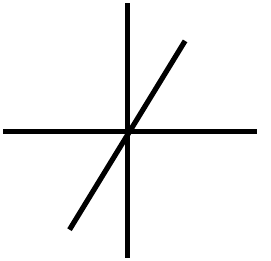}
		\put(-18.7,32){z}
		\put(-29,-1.5){x}
		\put(-.5,13){y}
	\end{array}$
\end{tabular}
\end{center}

First, rotate the piece of surface $180^\circ$ about the $z$-axis so that the `back' side is shown, then apply $d_c$.

\ 

{\footnotesize
\begin{tabular}{lll}
	$\begin{array}{c}
		\textrm{`front' side}\\
		\includegraphics[scale=1.3]{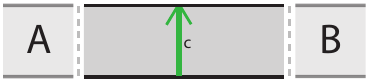}
	\end{array}$
& {\hskip -14 pt}
	$\begin{array}{c}
		\circlearrowright \textrm{about } z \textrm{ axis} \\
		-----  \rightarrow
	\end{array}$
& {\hskip -17 pt}
	$\begin{array}{c}
		\textrm{`back' side}\\
		\reflectbox{ \includegraphics[scale=1.3]{fbd1.pdf} }
	\end{array}$
\\
	\ 
& {\hskip -14 pt}
	$\begin{array}{c}
		\textrm{applying } d_c \\
		-----  \rightarrow
	\end{array}$
& {\hskip -14 pt}
	$\begin{array}{c}
		\ \\
		\includegraphics[scale=1.3]{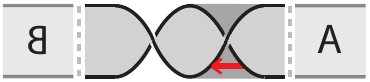}
	\end{array}$
\end{tabular}
}

\ 

Next we isotopically stretch \& twist parts of the surface so that what follows becomes easier to visualize.  What follows is a $180^\circ$ rotation about the $y$-axis that rotates everything except for the piece of surface marked with a *.

{\footnotesize
\begin{tabular}{lll}
	$\begin{array}{c}
		\ \\
		\includegraphics[scale=1.3]{fbd3.pdf}
	\end{array}$
& {\hskip -14 pt}
	$\begin{array}{c}
		\textrm{isotopy} \\
		----- \rightarrow
	\end{array}$
& {\hskip -14 pt}
	$\begin{array}{c}
		\ \\
		\includegraphics[scale=1.3]{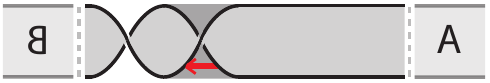}
	\end{array}$
\\
	\ 
& {\hskip -14 pt}
	$\begin{array}{c}
		\textrm{isotopy} \\
		----- \rightarrow
	\end{array}$
& {\hskip -14 pt}
	$\begin{array}{c}
		\ \\
		\includegraphics[scale=1.3]{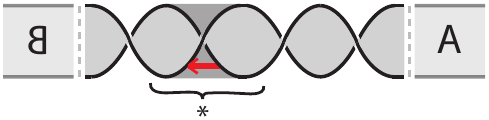}
	\end{array}$
\\
	\ 
& {\hskip -14 pt}
	$\begin{array}{c}
 		\circlearrowright \textrm{about } y \textrm{ axis} \\
		----- \rightarrow
	\end{array}$
& {\hskip -14 pt}
	$\begin{array}{c}
		\includegraphics[scale=1.3]{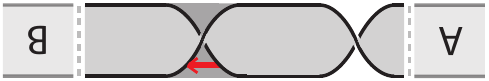}
	\end{array}$
\end{tabular}
}

\ 

Lastly, perform a $180^\circ$ rotation about the $x$-axis, then isotopically push down the crosscut to show that we have the same piece of surface that applying $d_c$ to the `front' side would have given.

{\footnotesize
{\hskip -14 pt}
\begin{tabular}{lll}
	$\begin{array}{c}
		\ \\
		\includegraphics[scale=1.3]{fbd6.pdf}
	\end{array}$  
& {\hskip -14 pt}
	$\begin{array}{c}
		\circlearrowright \textrm{about } x \textrm{ axis} \\
		----- \rightarrow
	\end{array}$  
& {\hskip -14 pt}
	$\begin{array}{c}
		\rotatebox{180}{\includegraphics[scale=1.3]{fbd6.pdf}}
	\end{array}$ 
\\
	\
& {\hskip -14 pt}
	$\begin{array}{c}
		\textrm{isotopy} \\
		-----  \rightarrow
	\end{array}$ 
& {\hskip -14 pt}
	$\begin{array}{c}
		\ \\
	\includegraphics[trim = 0mm 0mm 29mm 0mm, clip,scale=1.3]{fbd1.pdf}  \includegraphics[trim = 8.6mm 0mm 8.6mm 0mm, clip,scale=1.3]{fbd3.pdf}  \includegraphics[trim = 29mm 0mm 0mm 0mm, clip,scale=1.3]{fbd1.pdf}
	\end{array}$
\end{tabular} \qed
}
\\

We define the map $d$ in terms of the $d_c$'s and show that $d$ is well defined on $\mathcal{C}_{i,j,k,b}$.

\begin{defn}
\label{4.3}
Define $d:\mathcal{C}_{i,j,k,b} \longrightarrow \mathcal{C}_{i+1,j,k,b}$ by 

$$
d(F) := \sum_{\begin{array}{c} \textrm{\textcolor{darkgreen}{active} crosscuts} \\ c \in F \end{array}} (-1)^{\alpha(c)} d_c(F),
$$
where $\alpha(c)$ is the number of \textcolor{red}{inactive} crosscuts that come before $c$ in the ordering of crosscuts on $F$.  The map $d$ is defined on linear combinations of surfaces by linear extension.
\end{defn}

To see that the map $d$ is well defined, we must show that if $F$ is a $D^k$-surface with grading $(I,J,K,B) = (i,j,k,b)$, then $d(F)$ is a $D^k$-surface with grading $(I,J,K,B) = (i+1,j,k,b)$.

Before proving the well-definedness of $d$, a lemma concerning the generators for the $1^{\textrm{st}}$ homology of a compact surface will be useful.

\begin{lemma}
\label{4.4}
Given a compact surface $F$, and a neighborhood of a closed interval embedded inside the surface (as shown in the proof below), there is a basis for $H_1(F,\Q)$ which has at most one S.C.C. (simple closed curve) class representative running through the local piece of surface in question.  Furthermore, the basis can be chosen so that this representative only runs through this piece of surface once.
\end{lemma}

\emph{Proof.}  Let $F$ be a compact surface, suppose $\{ [\alpha_1], [\alpha_2], ... ,[\alpha_n] \}$ is a basis for $H_1(F)$ which has a S.C.C. class representative which runs through the piece of surface in question more than once.  Suppose without loss of generality that $\alpha_1$ is the S.C.C. class representative which runs through the piece of surface more than once.  Since $\alpha_1$ runs through the piece of surface multiple times, it may or may not change directions from one pass through the piece of surface to the next.

First we consider the case where $\alpha_1$ runs through the piece of surface at least once in both directions.
	\begin{center} 
	\includegraphics[trim = 0mm 0mm 0mm 28.8mm, clip,scale=1.3]{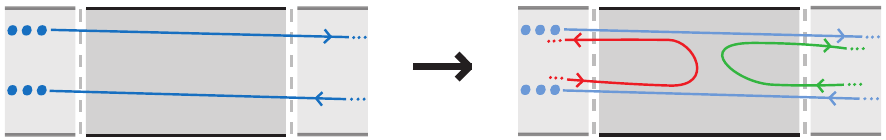}
	\put(-244,10){$\alpha_1$}
	\put(-244,43){$\alpha_1$}
	\put(-51,10){$\alpha_1$}
	\put(-51,43){$\alpha_1$}
	\put(-83,26){$\beta_0$}
	\put(-52,27){$\beta_1$}
	\end{center}
Replace $\alpha_1$ with the two S.C.C.'s shown above, $\beta_0$ and $\beta_1$, which are defined in terms of (part of) $\alpha_1$.  Construct $\beta_1$ and $\beta_0$ so that $[\alpha_1] = [\beta_1] - [\beta_0]$.  Then
$$
\textrm{span}(\{ [\alpha_1], [\alpha_2], ... ,[\alpha_n] \} ) \subseteq \textrm{span}(\{ [\beta_0], [\beta_1], [\alpha_2], ... , [\alpha_n]\} ).
$$
As $\{ [\alpha_1], [\alpha_2], ..., [\alpha_n] \}$ is a basis, it spans $H_1(F)$, so the above set inclusion is actually an equality.  Since $\{ [\beta_0], [\beta_1], [\alpha_2], ..., [\alpha_n]\}$ is a spanning set with $n+1$ elements, it must be linearly dependent.  This means that there is an $n$-element subset of $\{ [\beta_0], [\beta_1], [\alpha_2], ... ,[\alpha_n]\}$ which is a basis for $H_1(F)$.  Notice that $\beta_0$ and $\beta_1$ must run through the piece of surface strictly fewer times than $\alpha_1$ does.  Hence, this process can be repeated a finite number of times until there is a basis for $H_1(F)$ which contains no S.C.C.'s running through the piece of surface multiple times in \textit{different} directions.. 

Now consider the case where $\alpha_1$ runs through the piece of surface in the same direction each time.
	\begin{center} 
	\includegraphics[trim = 0mm 14.5mm 0mm 14.5mm, clip,scale=1.3]{gen321.pdf}
	\put(-244,12){$\alpha_1$}
	\put(-244,43){$\alpha_1$}
	\put(-50,12){$\alpha_1$}
	\put(-51,43){$\alpha_1$}
	\put(-100,19){$\uparrow$}
	\put(-101,8){$\beta_1$}
	\put(-67,18){$\uparrow$}
	\put(-68,7){$\beta_0$}
	\end{center}
Replace $\alpha_1$ with the two S.C.C.'s shown above, $\beta_0$ and $\beta_1$, which are defined in terms of (part of) $\alpha_1$.  This time, construct $\beta_1$ and $\beta_0$ so that $[\alpha_1] = [\beta_1] + [\beta_0]$, again yielding
$$
\textrm{span}(\{ [\alpha_1], [\alpha_2], ..., [\alpha_n] \} ) \subseteq \textrm{span}(\{ [\beta_0], [\beta_1], [\alpha_2], ... , [\alpha_n]\} ).
$$
As $\{ [\alpha_1], [\alpha_2], ..., [\alpha_n] \}$ is a basis for $H_1(F)$, it spans $H_1(F)$.  This implies that the above set inclusion is instead an equality.  The set $\{ [\beta_0], [\beta_1], [\alpha_2], ... , [\alpha_n]\}$ is then a spanning set for $H_1(F)$ with $n+1$ elements, and so it's linearly dependent.  By this linear dependence there must be an $n$-element subset which is a basis for $H_1(F)$.  Therefore, this process uses the basis $\{ [\alpha_1], [\alpha_2], ..., [\alpha_n] \}$ to find a basis that replaces $\alpha_1$ with a S.C.C. that passes through the piece of surface strictly fewer times than $\alpha_1$ does.  Repeat this process until a basis for $H_1(F)$ is found in which every S.C.C. passes through the given piece of surface at most once each.

Lastly, the following must be shown:  given a basis for $H_1(F)$ in which every S.C.C. passes through the given piece of surface 0 or 1 times, a basis can be created in which (at most) one S.C.C. runs through the piece of surface.  Let $\{ [\alpha_1], [\alpha_2], ..., [\alpha_n] \}$ be the given basis.  If 0 or 1 S.C.C.'s pass though the piece of surface, we are done.  Otherwise, suppose without loss of generality that $\alpha_1$ and $\alpha_2$ both pass through the piece of surface.  Whether or not $\alpha_1$ and $\alpha_2$ pass through the surface in the same direction does not matter, since $[\alpha_1] = [- \alpha_1]$.  
	\begin{center} 
	\includegraphics[trim = 0mm 28.8mm 0mm 0mm, clip,scale=1.3]{gen321.pdf}
	\put(-244,10){$\alpha_2$}
	\put(-244,43){$\alpha_1$}
	\put(-51,10){$\alpha_2$}
	\put(-51,43){$\alpha_1$}
	\put(-83,26){$\beta$}
	\put(-52,27){$\beta$}
	\end{center}
Now define $\beta$ as shown, so that $[\alpha_1] + [\alpha_2]= [\beta]$.  Since $[\alpha_1] = [\beta] - [\alpha_2]$, we have that 
$$
\textrm{span}(\{ [\alpha_1], [\alpha_2], ... , [\alpha_n] \} ) \subseteq \textrm{span}(\{ [\beta], [\alpha_2], ..., [\alpha_n]\} ).
$$
Since $\{ [\alpha_1], [\alpha_2], ..., [\alpha_n] \}$ is a basis for $H_1(F)$, the above set inclusion is actually an equality.  Thus, the set $\{ [\beta],[\alpha_2], ... , [\alpha_n] \}$ is a basis for $H_1(F)$ which has one fewer element passing through the given piece of surface than the basis $\{ [\alpha_1], [\alpha_2], ... , [\alpha_n] \}$.  This process replaces a basis which has two or more S.C.C.'s which pass through the piece of surface with a basis that has one fewer element passing through the piece of surface.  Repeat this process until we are left with a basis that has only one element passing through the given piece of surface.  \qed

Now we return to the matter at hand, working to show that $d$ is well defined. 

\begin{prop}
\label{4.5}
The differential $d$ increases the homological grading $I$ (= signature) by +1.
\end{prop}

\emph{Proof.}  It suffices to consider an arbitrary \textcolor{darkgreen}{active} crossing of an arbitrary surface $F \in \mathcal{C}_{i,j,k,b}$.
	\begin{center} 
	\includegraphics[scale=1.5]{dmap.pdf}
	\put(-228,14){$F$}
	\put(3,14){$d_c(F)$}
	\put(-27.5,.5){{\footnotesize \textbf{c} }}
	\end{center}
By Lemma \ref{4.4}, we may assume we have a basis $\{ [\beta_1], ... , [\beta_n] \}$ for $H_1(F, \Q)$ which has at most one curve, say $\beta_1$, passing through the crosscut $c$, at most once.  Since $F$ and $d_c(F)$ only differ locally at the piece of surface shown, we can construct a basis $\{ [\alpha_0], [\alpha_1], ... , [\alpha_n] \}$ for $H_1(d_c(F), \Q)$ from the basis $\{ [\beta_1], ... , [\beta_n] \}$.  Since none of $\beta_2, ... ,\beta_n$ cross the local piece of surface, we can choose $\{ \alpha_2, ..., \alpha_n \}$ that are isotopic to $\{ \beta_2, ..., \beta_n \}$.  Let $\alpha_1$ be the curve corresponding to $\beta_1$ as shown below.
	\begin{center} 
	\includegraphics[scale=1.5]{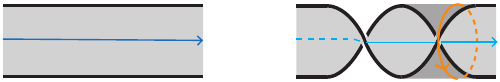}
	\put(-228,14){$F$}
	\put(-177,22){$\beta_1$}
	\put(3,14){$d_c(F)$}
	\put(-49,20){$\alpha_1$}
	\put(-33,0){$\alpha_0$}
	\end{center}
Notice that a curve traveling around the $\frac{1}{2}$-twisted loop in $d_c(F)$ is needed to span $H_1(F, \Q)$.  Let $\alpha_0$ be the curve traveling around the $\frac{1}{2}$-twisted loop shown in the figure above.  Then $\{ [\alpha_0], [\alpha_1],..., [\alpha_n] \}$ is a basis for $H_1(d_c(F), \Q)$.

Now consider the Goeritz matrices for $d_c(F)$ and $F$.  We have that dim$(G_{d_c(F)}) = (n+1) \times (n+1)$ and dim$(G_F) = n \times n$.  Let $m = \textrm{lk} (\alpha_1, \tau \alpha_1)$.  Since $\alpha_1$ travels along surface that has an extra right-handed $\frac{1}{2}$-twist compared with $\beta_1$, we have that lk$(\beta_1, \tau \beta_1) = m-1$.  Let $A$ be the $(n-1) \times (n-1)$ minor of the Goeritz matrix corresponding to $\{ \alpha_2, ..., \alpha_n \}$.  Since $\{ \alpha_2, ..., \alpha_n \}$ are isotopic to $\{ \beta_2, ..., \beta_n \}$, $A$ is also the minor of the Goeritz matrix corresponding to $\{ \beta_2, ..., \beta_n \}$.  From this we have that
$$
\left( \ G_{d_c(F)} {\hskip -12 pt} \ ^{\ ^{\ ^{\ }}}\right) 
\ = \ \ \ \ \ 
\left( \begin{array}{c|c|c}
1&1&0 \\
\hline
1&m&\begin{array}{lll}
\leftarrow&\hskip -8 pt R \hskip -8 pt&\rightarrow
\end{array} \\
\hline
0&\begin{array}{c}
\uparrow \\
C \\
\downarrow
\end{array}&A
\end{array} \right)
	\put(-123,26){$\alpha_0$}
	\put(-123,13){$\alpha_1$}
	\put(-123,0){$\alpha_2$}
	\put(-120,-14){$\vdots$}
	\put(-123,-23){$\alpha_n$}
	\put(-107,38){$\tau \alpha_0$}
	\put(-83,38){$\tau \alpha_1$}
	\put(-57,38){$\tau \alpha_2$}
	\put(-39,38){$\cdots$}
	\put(-25,38){$\tau \alpha_n$}
\ \ \textrm{ and } \ \ 
\left( \ G_F {\hskip -12 pt} \ ^{\ ^{\ ^{\ }}} \right) 
\ = \ \ \ \ \ 
\left( \begin{array}{c|c}
m-1&\begin{array}{lll}
\leftarrow&\hskip -8 pt R \hskip -8 pt&\rightarrow
\end{array} \\
\hline
\begin{array}{c}
\uparrow \\
C \\
\downarrow
\end{array}&A
\end{array} \right)
	\put(-116,18){$\beta_1$}
	\put(-116,5){$\beta_2$}
	\put(-113,-9){$\vdots$}
	\put(-116,-19){$\beta_n$}
	\put(-87,31){$\tau \beta_1$}
	\put(-56,31){$\tau \beta_2$}
	\put(-40,31){$\cdots$}
	\put(-26,31){$\tau \beta_n$},
$$
where row and column of entries denoted by $R$ and $C$ are the same in the two matrices above.  This is due to the fact that $\alpha_1$ and $\beta_1$ are the same away the local pieces of surface shown.  Now, by Sylvester's law of inertia \cite{Sy}, signature is unchanged by matrix congruence, and so the following calculation is useful (where `$\cong$' denotes matrix congruency, not matrix similarity):
$$
\left( \ G_{d_c(F)} {\hskip -12 pt} \ ^{\ ^{\ ^{\ }}}\right) 
=
\left( \begin{array}{c|c|c}
1&1&0 \\
\hline
1&m&\begin{array}{lll}
\leftarrow&\hskip -8 pt R \hskip -8 pt&\rightarrow
\end{array} \\
\hline
0&\begin{array}{c}
\uparrow \\
C \\
\downarrow
\end{array}&A
\end{array} \right)
\cong
\left( \begin{array}{c|c|c}
1&0&0 \\
\hline
0&m-1&\begin{array}{lll}
\leftarrow&\hskip -8 pt R \hskip -8 pt&\rightarrow
\end{array} \\
\hline
0&\begin{array}{c}
\uparrow \\
C \\
\downarrow
\end{array}&A
\end{array} \right)
=
\left( \begin{array}{c|c}
1&0 \\
\hline
0&G_F
\end{array} \right).
$$

It follows from the divide-and-conquer method for computing eigenvalues that sig($d_c(F)$) = sig($F$)+1. \qed

\begin{prop}
The differential $d$ fixes the gradings $J,K$ and $B$.
\end{prop}

\emph{Proof.}  It is a simple exercise to check that applying $d_c$ changes the Euler characteristic by -1.  Since applying $d_c$ does not change the number of dots on a surface, we have that $\delta(F) = \delta(d_c(F))$.  Proposition \ref{4.5} implies that applying $d_c$ increases signature by +1, and so $J =- \chi - I + 2\delta$ the same for $F$ and $d_c(F)$.  

The number of crosscuts does not change, so $K$ is also fixed.  

Finally, using Proposition \ref{4.5} gives us that $B = (b$ + \# of \textcolor{darkgreen}{active} crosscuts) is constant. \qed

\begin{prop}
The differential $d$ applied to a $D^k$-surface is a $D^k$-surface.
\end{prop}

\emph{Proof.}  Assume $F$ is a $D^k$-surface, with oriented crosscuts $\{c_1, ...,c_k\}$.  Consider the cross-dual surfaces of $F$ and $d_c(F)$ as well as their skeletons.  Since $F$ and $d_c(F)$ only differ near the crosscut $c$, and since the corresponding cross-dual surfaces skeletons only potentially differ in that same area, we restrict our attention to a neighborhood of the crosscut $c$.
	\begin{center} 
	\includegraphics[scale=1.5]{dmap.pdf} \put(-228,14){$F$} \put(3,14){$d_c(F)$} \\
	\vspace{3 pt}
	$\downarrow \textrm{cross-dual} \downarrow \ \ \ \ \ \ \ \ \ \ \ \ \ \ \ \ \ \ \ \ \downarrow\textrm{cross-dual} \downarrow$\\
	\vspace{7 pt}
	\includegraphics[scale=1.5]{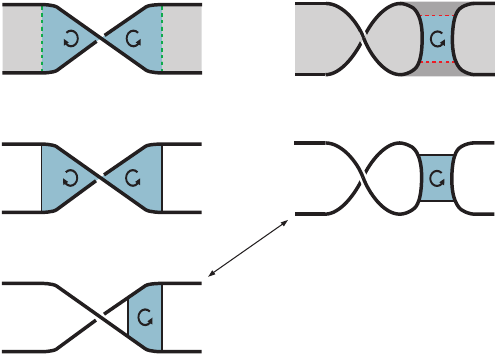}
	\put(-233,133){$F^{\textrm{cd}}$}
	\put(2,133){$(d_c(F))^{\textrm{cd}}$}
	\put(-196.5,103.5){$\downarrow \textrm{skeleton} \downarrow \ \ \ \ \ \ \ \ \ \ \ \ \ \ \ \ \ \ \ \ \ \ \  \downarrow\textrm{skeleton} \downarrow$}
	\put(-194.5,43){$\downarrow \textrm{isotopy} \downarrow$}
	\put(-122,29){ \rotatebox{35}{isotopic} }
	\\
	\end{center}
Consider the cross-dual surfaces shown above.  If we forget the local orientations then label the two sides of the surface coming in from the left and the right of each figure, we find that $\fcd$ is orientable if and only if $(d_c(F))^{\textrm{cd}}$ is orientable.  Also, note that the skeletons of $\fcd$ and $(d_c(F))^{\textrm{cd}}$ are isotopic to one another (even with the consideration of the local orientations).  Therefore, if $F$ satisfies the conditions required to be a $D^k$-surface, then $d_c(F)$ satisfies those conditions as well.
\qed

\bigskip
The previous propositions together prove the following theorem.

\begin{theorem}
The map $d:\mathcal{C}_{i,j,k,b} \longrightarrow \mathcal{C}_{i+1,j,k,b}$ is well-defined.
\end{theorem}

Now recall each local map $d_c$ acts on a neighborhood of the crosscut $c$ and nowhere else.  Given two crosscuts $c_i$ and $c_j$, we can choose small enough neighborhoods so that the maps $d_{c_i}$ and $d_{c_j}$ do not interact with one another.  This means that these maps (positively) commute with one another.  Since we defined the map $d:\mathcal{C}_{i,j,k,b} \longrightarrow \mathcal{C}_{i+1,j,k,b}$ by 
$$
d := \sum_{\begin{array}{c} F \in \mathcal{C}_{i,j,k,b} \end{array} } \sum_{\begin{array}{c} \textrm{\textcolor{darkgreen}{active} crosscuts} \\ c \in F \end{array}} (-1)^{\alpha(c)} d_c,
$$
we get that all pairs of maps $d_{c_i}$ and $d_{c_j}$ will negatively commute.  This is because $(-1)^{\alpha(c_j)}d_{c_j}\circ (-1)^{\alpha(c_i)}d_{c_i}$ and $(-1)^{\alpha(c_i)}d_{c_i}\circ (-1)^{\alpha(c_j)}d_{c_j}$ will have opposite sign regardless of whether $i<j$ or $j<i$.  Thus the following theorem is proved.

\begin{theorem}
$d \circ d = 0.$
\end{theorem}

\section{Reduced Chain Complexes}

In this section, we use the chain complex defined in the previous sections to define a reduced chain complex.  The idea is to take the quotient of the chain complex $\mathcal{C}$ by a certain submodule.  This submodule is generated by skein relations that come from a Frobenius system.

Before defining Frobenius system, we must first define Frobenius extension.  There are different (equivalent) ways to define Frobenius extension, but we will follow \cite{Kh2}, where Khovanov defines a \textit{Frobenius extension} as an inclusion $\iota :R \rightarrow A$ of commutative unital rings whose left and right adjoint functors are isomorphic.  Khovanov continues by giving the following proposition (whose proof can be found in section 4 of \cite{Kad}).

\begin{prop}
The inclusion $\iota$ is a Frobenius extension if and only if there exists an $A$-bimodule map $\Delta : A \rightarrow A \otimes_R A$ and an $R$-module map $\varepsilon : A \rightarrow R$ such that $\Delta$ is coassociative and cocommutative, and $(\varepsilon \otimes \textrm{Id})\Delta = \textrm{Id}$.
\end{prop}

Still following \cite{Kh2}, we now define Frobenius system with the above proposition in mind.

\begin{defn}
A Frobenius extension, together with a choice of $\varepsilon$ and $\Delta$, will be denoted $\mathcal{F} = (R,A,\varepsilon,\Delta)$ and called a \textit{Frobenius system}.
\end{defn}

\begin{remark}
Since we will be using results from Kaiser's paper \cite{Kai}, it should be noted that Kaiser defines a Frobenius algebra to be what we call a Frobenius system.  Kaiser reserves the term Frobenius system for when a choice of $u_i,v_i \in A$ for which $\Delta ( \1 ) = \Sigma_{i=1}^r u_i \otimes_R v_i$ is specified.
\end{remark}

Now we turn our attention to a specific Frobenius system (which we will denote by $\mathcal{F}_5$ in order to be consistent with \cite{Kh2}).

\begin{defn}
Let $\mathcal{F}_5$ be the Frobenius system defined by $R = \Z[h,t], A = R[x]/(x^2 - hx - t)$, $\varepsilon (\1) = 0, \ \varepsilon (x) = \1$, and $\Delta(\1) = \1 \otimes x + x \otimes \1 - h (\1 \otimes \1)$.  Then $\Delta(x)$ follows from $A$-bilinearity:
\begin{center}
$ \begin{array}{rl}
\vspace{2 pt} \Delta(x) \ = \ \Delta(\1) x \hspace{-4 pt} & = \ (\1 \otimes x + x \otimes \1 - h(\1 \otimes \1) ) x\\
 \vspace{2 pt} & = \ \1 \otimes x^2 + x \otimes x - h(\1 \otimes x )\\ 
 \vspace{2 pt} & = \ \1 \otimes (hx + t) + x \otimes x - \1 \otimes (hx)\\
& = \ \1 \otimes t + x \otimes x 
\end{array} $
\end{center}
\end{defn}

\begin{remark}
The Frobenius system $\mathcal{F}_5$ is referred to as the the \textit{rank two universal} Frobenius algebra in \cite{Kai}.  It is universal in the sense that any rank two Frobenius system can be obtained from $\mathcal{F}_5$ by base changing and `twisting' -- see \cite{Kh2} for further explanation and a proof.
\end{remark}

In \cite{Kai}, Kaiser obtains skein relations from a given Frobenius system by allowing surfaces to be `colored' by elements of $A$.  We follow Kaiser's procedure using $\mathcal{F}_5$ as our Frobenius system.  Since $\{ \1, x \}$ is a basis for $\mathcal{F}_5$, it suffices to only consider colorings by $\1$ and $x$.  In this paper, $D^k$-surfaces may be marked by dots.  Identify (the coloring of) $\1$ to all undotted facets of $D^k$-surfaces and identify (the coloring of) $x$ to all once-dotted facets of $D^k$-surfaces.  

The $D^k$-surfaces with $n \geq 2$ dots would be identified with a coloring of $x^n$; however, since $x^2-hx-t=0$, such surfaces may always be written as $\Z[h,t]$-linear combinations of surfaces with at most one dot on each facet.  An example of this is shown below.

\begin{example}
\label{d2}
$$
		\includegraphics[scale=1.5]{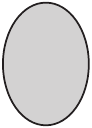} 
		\thicklines
		\color{red}
		\put(-37.5,22){\vector(3,1){35.3}}
		\color{blue}
		\put(-28,40){$\bullet$}
		\put(-18,38){$\bullet$}
		\put(-22,13){$\bullet$}
		\color{black}
\ \ 
	\begin{array}{c}
		= \\ \ \\ \ \\ \ \\ \ 
	\end{array}
	\begin{array}{c}
		h\\ \ \\ \ \\ \ \\ \ 
	\end{array}
		\includegraphics[scale=1.5]{u0.pdf} 
		\thicklines
		\color{red}
		\put(-37.5,22){\vector(3,1){35.3}}
		\color{blue}
		\put(-18,38){$\bullet$}
		\put(-22,13){$\bullet$}
		\color{black}
\ \ 
	\begin{array}{c}
		+ \\ \ \\ \ \\ \ \\ \ 
	\end{array}
	\begin{array}{c}
		t\\ \ \\ \ \\ \ \\ \ 
	\end{array}
		\includegraphics[scale=1.5]{u0.pdf} 
		\thicklines
		\color{red}
		\put(-37.5,22){\vector(3,1){35.3}}
		\color{blue}
		\put(-22,13){$\bullet$}
		\color{black}
$$
\end{example}
\vspace{-40 pt}

\begin{remark}
By assigning one of the two basis elements $\{ \1, x \}$ of $A = \Z[x]/(x^2-hx-t)$ to every facet of each (non-zero) $D^k$-surface $F$, we are associating $A^{\otimes n}$ to each $D^k$-surface $F$ with $n$ facets.  This will be useful later when calculating homology.  We have:  
$$
\mathcal{C}_{i,j,k,b}(L) := \bigoplus_{ D^k\textrm{-surfaces }F\textrm{, with} } A^{\otimes n},
\put(-108,-22){ {\scriptsize $I = i, J = j, B= b,\partial(F) = L$} }
$$
where $n$ is the number of facets of $F$.
\end{remark}

Following \cite{Kai}, we now use our Frobenius system to define the following submodules.

\begin{defn}
\label{relations}
Let $\mathfrak{R}(\mathcal{F}_5)_{i,j,k,b}(L)$ be the submodule of $\mathcal{C}_{i,j,k,b}(L)$ generated by the following four elements:
\begin{itemize}
\item [(S0)] Any $D^k$-surface $F \in \mathcal{C}_{i,j,k,b}(L)$ which has an undotted sphere as a component.  That is, there exists a (possibly null) $D^k$-surface $F' \in \mathcal{C}_{i,j,k,b}(L)$ such that $F = F' \sqcup \Sigma$, where $\sqcup$ denotes disjoint union and $\Sigma$ is an undotted sphere.
\begin{center}
\begin{tabular}{cc}
$\ F' \ \ \ \ \sqcup$
&
	$\begin{array}{l}
		\includegraphics[trim = 27.5mm 0mm 0mm 0mm, clip,scale=.38]{xyz_sph.pdf}
	\end{array}$
\end{tabular}
\put(-11,-2){ {\Huge )} }
\put(-100.5,-2){ {\Huge (} }
\end{center}
\item [(S1)] The difference of a $D^k$-surface $F \in \mathcal{C}_{i,j,k,b}(L)$ and the union of $F$ with a once-dotted sphere component.
\begin{center}
\begin{tabular}{ccc}
$F \ \ \ \ \ \ - \ \ $
&
$F \ \ \ \ \sqcup$
&
	$\begin{array}{l}
		\includegraphics[trim = 27mm 0mm 0mm 0mm, clip,scale=.4]{xyz_sph.pdf}
		\put(-17,21){$\textcolor{blue}{\bullet}$}
	\end{array}$
\end{tabular}
\put(-11,-2.5){ {\Huge )} }
\put(-100.5,-2.5){ {\Huge (} }
\put(-133,-2.5){ {\Huge )} }
\put(-155.5,-2.5){ {\Huge (} }
\end{center}
\item [(NC)] If a $D^k$-surface $F \in \mathcal{C}_{i,j,k,b}(L)$ has a simple closed curve $\gamma \in F$ that $[1]$ does not intersect any crosscuts and $[2]$ bounds a disk $D \in S^3$ satisfying $D \cap F = \gamma$, then let $\bar{F}$ denote the surface obtained from compressing $F$ along $D$.  Compression involves replacing an annular neighborhood of $\gamma$ with two disks, $D_-$ and $D^+$, a process that may or may not split the facet involved into two different facets.  Let $\bar{F}_\bullet$ (resp. $\bar{F}^\bullet$) denote the surface $\bar{F}$ with a dot placed on the part of the surface comprised of $D_-$ (resp. $D^+$).  The combination $F - \bar{F}_\bullet - \bar{F}^\bullet - h \bar{F}$ is an element of $\mathfrak{R}(\mathcal{F}_5)_{i,j,k,b}(L)$. 

A local picture of such an element is given below.
	\begin{center} 
	\includegraphics[trim = 0mm 0mm 10.5mm 0mm, clip,scale=1.3]{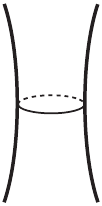}
	\begin{tabular}{c}
		$ - $ \\ \ \\ \ \\ \ \\ \ \\ \ \\ \ 
	\end{tabular}
	\includegraphics[trim = 10.5mm 0mm 0mm 0mm, clip,scale=1.3]{nc2.pdf}
	\put(-21,28){$\textcolor{blue}{\bullet}$}
	\begin{tabular}{c}
		$ - $ \\ \ \\ \ \\ \ \\ \ \\ \ \\ \ 
	\end{tabular}
	\includegraphics[trim = 10.5mm 0mm 0mm 0mm, clip,scale=1.3]{nc2.pdf}
	\put(-21,43){$\textcolor{blue}{\bullet}$}
	\begin{tabular}{c}
		$ - \ \ \ h$ \\ \ \\ \ \\ \ \\ \ \\ \ \\ \ 
	\end{tabular}
	\hspace{-12 pt}\includegraphics[trim = 10.5mm 0mm 0mm 0mm, clip,scale=1.3]{nc2.pdf}
	\end{center}
	\vspace{-45 pt}
\end{itemize}
\end{defn}

\begin{remark}
Since we plan to take the quotient $\mathcal{C}_{i,j,k,b}(L) / \mathfrak{R}(\mathcal{F}_5)_{i,j,k,b}(L)$, the elements of $\mathfrak{R}(\mathcal{F}_5)_{i,j,k,b}(L)$ will be equal to zero.  Hence, we will often refer to (S0), (S1) and (NC) as `relations' (specifically, they are known as the sphere relations and the neck-cutting relation, respectively).  The reader should note that for relations involving multiple $D^k$-surfaces, if one surface is a $D^k$-surface, then all surfaces are $D^k$-surfaces.  To ensure that $\mathfrak{R}(\mathcal{F}_5)_{i,j,k,b}(L)$ is indeed a submodule of $\mathcal{C}_{i,j,k,b}(L)$, we must have that relations involving multiple $D^k$-surfaces all have the same values of the gradings $I,J,K$ and $B$.  This is proved below.
\end{remark}

\begin{prop}
For each relation defined in $\mathfrak{R}(\mathcal{F}_5)_{i,j,k,b}(L)$, all surfaces involved have the same value of $I, J, K, $ and $B$. 
\end{prop}

\emph{Proof.}  It is only necessary to check the relations involving multiple surfaces, (S1) and (NC).
\begin{itemize}
\item [(S1):] A sphere cannot contribute to signature, so $I$ is unchanged.  Since a dotted sphere contributes +2 to Euler characteristic and has one dot, $J =- \chi - I + 2\delta$ is the same after the removal of the dotted sphere.  As no crosscuts are present on a sphere, $K$ and $B$ are also unchanged.
\item [(NC):] For the index $I$, consider the generators for the $1^{\textrm{st}}$ homology of the surfaces involved.  If the surface with the neck does not need to have a generator run along the neck, then a curve traveling around the neck must be null-homotopic.  In this case the local relation would not affect any of the generating curves, and hence would not affect signature.  Suppose there are generators of the $1^{\textrm{st}}$ homology that run along the neck.  By an argument similar to that in the proof of Lemma \ref{4.4}, we can find a basis $\{ [\alpha_1], [\alpha_2], ... , [\alpha_n] \}$ in which only one element runs along the neck, traveling across it only once.  Without loss of generality, call this element $\alpha_1$.  In this case, a curve traveling around the neck could \textit{not} be null-homotopic, and therefore would be a non-trivial element of the $1^{\textrm{st}}$ homology, different from $\alpha_1$.  Without loss of generality, let $\alpha_2$ be this curve.
	\begin{center} 
	\includegraphics[scale=2]{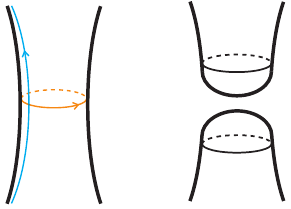}
	\put(-151,92){$\alpha_1$}
	\put(-131,51){$\alpha_2$}
	\end{center}
Now, since the set $\{ [\alpha_1], [\alpha_2],[\alpha_3], ... , [\alpha_n] \}$ is a basis for the $1^{\textrm{st}}$ homology of the surface with the neck, then $\{ [\alpha_3], ... , [\alpha_n] \}$ would serve as a basis for the $1^{\textrm{st}}$ homology of the corresponding surface with the cut neck.  Now let $G_{\textrm{neck}}$ be the Goeritz matrix for the surface with the neck, and $G_{\textrm{cut}}$ be the Goeritz matrix for the corresponding surface with the cut neck.  Then dim($G_{\textrm{neck}}$) = $n \times n$ and dim($G_{\textrm{cut}}$) = $(n-2) \times (n-2)$.  Letting $m = \textrm{lk} (\alpha_1, \tau \alpha_1)$, we have the following matrix congruence calculation:
$$
\left( \ G_{\textrm{neck}} {\hskip -12 pt} \ ^{\ ^{\ ^{\ }}}\right) 
\ = \ \ \ \ 
\left( \begin{array}{c|c|c}
m&1&\begin{array}{lll}
\leftarrow&\hskip -8 pt ? \hskip -8 pt&\rightarrow
\end{array} \\
\hline
1&0&0 \\
\hline
\begin{array}{c}
\uparrow \\
? \\
\downarrow
\end{array}&0&G_{\textrm{cut}}
\end{array} \right)
	\put(-118,25){$\alpha_1$}
	\put(-118,12.5){$\alpha_2$}
	\put(-118,0){$\alpha_3$}
	\put(-115,-14){$\vdots$}
	\put(-118,-23){$\alpha_n$}
	\put(-97,38){$\tau \alpha_1$}
	\put(-75,38){$\tau \alpha_2$}
	\put(-55,38){$\tau \alpha_3$}
	\put(-38,38){$\cdots$}
	\put(-25,38){$\tau \alpha_n$}
\cong
\left( \begin{array}{c|c|c}
0&1&0 \\
\hline
1&0&0 \\
\hline
0&0&G_{\textrm{cut}}
\end{array} \right)
\cong
\left( \begin{array}{c|c|c}
-1&0&0 \\
\hline
0&1&0 \\
\hline
0&0&G_{\textrm{cut}}
\end{array} \right).
$$
By Sylvester's law of inertia \cite{Sy} and the divide-and-conquer method for computing eigenvalues, it follows that sig$(G_{\textrm{neck}})$ = sig$(G_{\textrm{cut}})$. 

Next, it is straightforward to calculate that cutting a neck changes the Euler characteristic by +2.  Since the surface with the cut neck gets an additional dot, we have that $J =- \chi - I + 2\delta$ is unchanged.  Finally, the local relation $d$ does not affect crosscuts, so $K$ and $L$ are also unchanged.  \qed
\end{itemize}

\begin{defn} 
Define the {\it universal chain modules} by $\UC_{i,j,k,b}(L) := \mathcal{C}_{i,j,k,b}(L) / \mathfrak{R}(\mathcal{F}_5)_{i,j,k,b}(L)$.  Hence, each $\UC_{i,j,k,b}(L)$ is a free module of equivalence classes of isotopy classes of $D^k$-surfaces in $S^3$ with $I=i, J=j, K=k,B=b,$ and $\partial(F) = L$.  Let $\UC$ denote the resulting chain complex, called the \textit{universal chain complex}.
\end{defn}

The map $d$ will also serve as the differential for the universal chain complex $\UC$.  For $d$ to be well defined on $\UC$, we must have that $d(F) \in \mathfrak{R}(\mathcal{F}_5)_{i,j,k,b}(L)$, for all $F \in \mathfrak{R}(\mathcal{F}_5)_{i,j,k,b}(L)$.  It suffices to show that applying $d$ to each of the three types of relations that generate $\mathfrak{R}(\mathcal{F}_5)_{i,j,k,b}(L)$ is equal to zero.  This is proved below.

\begin{prop}
The differential $d$ applied to each relation equals zero.
\label{4.8}
\end{prop}

\emph{Proof.}  First consider the (S0) relation.  Note that $d_c$ only changes a surface near a neighborhood of a crosscut.  Since a sphere cannot have any crosscuts on it, a surface which has a sphere bounding a ball as a component will still have this same sphere after applying $d_c$.  Therefore the relation still equals zero after applying $d_c$.

For similar reasons, a sphere with a dot will not be affected by applying $d_c$.  Thus a dotted sphere may be removed before or after applying $d_c$, with the same effect.  So (S1) is also zero after applying $d_c$.

Finally, for (NC), we again exploit the fact that the all of the action of $d_c$ happens away from the relation in question.  Recall that for the (NC) relation, compressing along a disk whose boundary intersects a crosscut is not permitted.  Hence the equality of the surface with the neck and the sum of the dotted surfaces without the neck carries through after applying $d_c$.  \qed

\bigskip
It follows that we still have $d \circ d = 0$, so the universal chain complex is indeed a chain complex.  Though one could proceed using the full generality of the $\mathcal{F}_5$ (the `universal rank two' Frobenius system), we will only consider the special case of when $h=0$ and $t=0$, where $A = \Z[h,t][x]/(x^2-hx-t) = \Z[x]/(x^2)$.  This gives the Frobenius system introduced by Khovanov in \cite{Kh1}.

\begin{defn}
Let $\mathcal{F}_1$ be the Frobenius system $\mathcal{F}_5$, but with $h=0$ and $t=0$ (We use $\mathcal{F}_1$ to match the notation used in \cite{Kh2}).  That is, $\mathcal{F}_1$ is the Frobenius system with $R = \Z, A = \Z[x]/(x^2)$, $\varepsilon(\1)=0, \ \varepsilon(x)=\1, \Delta(\1) = \1 \otimes x + x \otimes \1$, and $\Delta(x) = x \otimes x$.
\end{defn}

The resulting relations that generate the submodule $\mathfrak{R}(\mathcal{F}_1)_{i,j,k,b}(L)$ of $\mathcal{C}$ are referred to as \textit{Bar-Natan skein relations} by Asaeda and Frohman in \cite{A-F}.  Hence, we will name the quotient modules that we obtain accordingly.

\begin{defn}
\label{complex_defn}
Define the {\it Bar-Natan chain modules} by $\BNC_{i,j,k,b}(L) := \mathcal{C}_{i,j,k,b}(L) / \mathfrak{R}(\mathcal{F}_1)_{i,j,k,b}(L)$.  Hence, each $\BNC_{i,j,k,b}(L)$ is a free module of equivalence classes of isotopy classes of $D^k$-surfaces in $S^3$ with $I=i, J=j, K=k,B=b,$ and $\partial(F) = L$.  Let $\BNC$ denote the resulting chain complex, called the \textit{Bar-Natan chain complex}.
\end{defn}

\begin{remark}
The Frobenius system $\mathcal{F}_1$ is the Frobenius system that will be used for the remainder of the paper.  Since using $\mathcal{F}_1$ implies that $h=0$ and $t=0$, the (NC) relation (from Definition \ref{relations}) and the process of reducing facets with 2 or more dots (as shown in Example \ref{d2}) are both simplified.  Below we give a visual summary of how the relations from Definition \ref{relations} and the dot-reducing process appear under the use of the Frobenius system $\mathcal{F}_1$.  From here on, we will refer to this dot-reducing process as the (D2) relation.  

Together, these four relations are known as the Bar-Natan skein relations in \cite{A-F}.
\begin{itemize}
\item [(S0)] \ 
\begin{center}
\begin{tabular}{cc}
$\ F \ \ \ \ \sqcup$
&
	$\begin{array}{l}
		\includegraphics[trim = 27.5mm 0mm 0mm 0mm, clip,scale=.38]{xyz_sph.pdf}
	\end{array}$ \ \ \ \ \ \ \ \ \ \ 
\end{tabular}
\put(-44,-2){ {\Huge )} }
\put(-134,-2){ {\Huge (} }
\put(-34,1.5){ \ \ = \ \ 0}
\end{center}
\item [(S1)] \ 
\begin{center}
\begin{tabular}{ccc}
$F \ \ \ \ \sqcup$
&
	$\begin{array}{l}
		\includegraphics[trim = 27mm 0mm 0mm 0mm, clip,scale=.4]{xyz_sph.pdf}
		\put(-17,21){$\textcolor{blue}{\bullet}$}
	\end{array}$
&
$\ = \ \ \ \ \ F$
\end{tabular}
\put(-7,-2.5){ {\Huge )} }
\put(-28,-2.5){ {\Huge (} }
\put(-60,-2.5){ {\Huge )} }
\put(-152.5,-2.5){ {\Huge (} }
\end{center}
\item [(NC)] \ 
	\begin{center} 
	\includegraphics[trim = 0mm 0mm 10.5mm 0mm, clip,scale=1.3]{nc2.pdf}
	\begin{tabular}{c}
		$ = $ \\ \ \\ \ \\ \ \\ \ \\ \ \\ \ 
	\end{tabular}
	\includegraphics[trim = 10.5mm 0mm 0mm 0mm, clip,scale=1.3]{nc2.pdf}
	\put(-21,28){$\textcolor{blue}{\bullet}$}
	\begin{tabular}{c}
		$ + $ \\ \ \\ \ \\ \ \\ \ \\ \ \\ \ 
	\end{tabular}
	\includegraphics[trim = 10.5mm 0mm 0mm 0mm, clip,scale=1.3]{nc2.pdf}
	\put(-21,43){$\textcolor{blue}{\bullet}$}
	\end{center}
	\vspace{-41 pt}
\item [(D2)] \ 
\begin{center}
	$\begin{array}{l}
		\includegraphics[scale=.4]{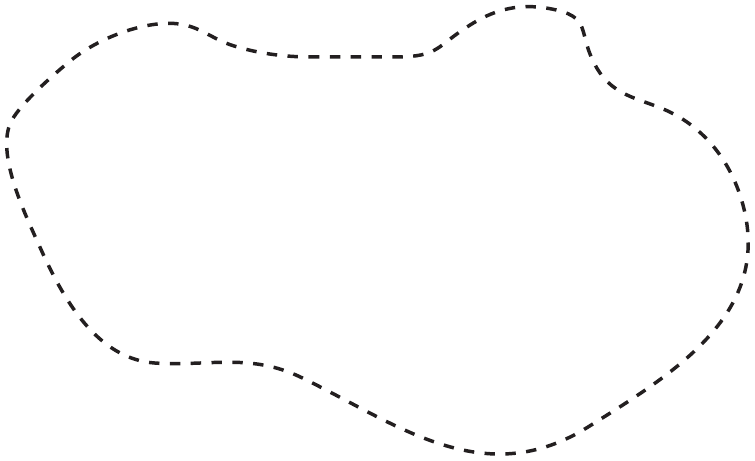}
	\end{array}$ \ \ \ \ \ \ 
    \put(-3,1.5){ \ \ = \ \ 0}
	\put(-58,8){$\textcolor{blue}{\bullet}$}
	\put(-38,0){$\textcolor{blue}{\bullet}$}
\end{center}
\end{itemize}
\end{remark}

\section{The Diagramless Homology of a Link \textit{L}}
\label{dhsection}

In Section \ref{2}, the crosscuts of a $D^k$-surface are required to be ordered.  In this section, we see that different choices of crosscut orderings yield the same homology for $L$ (see Theorem \ref{ccthm1}).  In fact, our Bar-Natan complex can often be seen as a sum of distinct (but isomorphic) subcomplexes, each of which corresponds to different crosscut ordering.  The resulting homology can usually be seen some number of repeat copies of a more basic quantity. 

We wish to avoid having multiple copies of the same homology due to different (but equivalent) crosscut orderings.  With this in mind, we will define an equivalence relation between certain subcomplexes which only differ by a crosscut reordering.  The \textit{diagramless homology} of a link $L$, denoted $D\mathcal{H}(L)$, will be constructed by choosing one representative from each equivalence class of subcomplexes for the diagramless complex.

\begin{theorem}
\label{ccthm1}
Let $\hspace{1pt} \BNC$ be a Bar-Natan chain complex for the link $L$ with a certain ordering of crosscuts for each $D^k$-surface (equivalence class) in $\BNC$, and let $\hspace{1pt} \BNC'$ be the Bar-Natan chain complex for the link $L$ differing from $\hspace{1pt} \BNC$ only in that the ordering of crosscuts has been permuted.  Then there exists a chain isomorphism $\psi: \BNC \rightarrow \BNC'$.
\end{theorem}

\emph{Proof.}  Without loss of generality, label the crosscuts of each $D^k$-surface (equivalence class) of $\BNC$ by $\{ c_1, c_2, ... , c_k \}$.  The complex $\BNC'$ differs from $\BNC$ only in that the ordering of the crosscuts has been permuted.  This reordering of the crosscuts can be seen as an action of the symmetric group $S_k$ on the subscripts of the crosscut labels.  Given $\sigma \in S_k$, we have $\sigma( \{ c_1, c_2, ... , c_k \} ) = \{ c_{\sigma(1)}, c_{\sigma(2)}, ... , c_{\sigma(k)}  \}$.  Since $S_k$ is generated by all the transpositions of the form ($i$ $i$+1), for $i=1,...,k-1$, it suffices to prove that the reordering given by  $\sigma =$ ($i$ $i$+1) yields an isomorphic chain complex.  

For each crosscut, let $0$ denote an \textcolor{darkgreen}{active} crosscut and let $1$ denote an \textcolor{red}{inactive} crosscut.  Let $m_1 m_2 \cdots m_k$ denote a length $k$ string of $0$'s and $1$'s -- we'll call this a \textit{binary string}.  To each surface in $\BNC$ we associate the binary string $m_1 m_2 \cdots m_k$, where $m_j = 0$ if $c_j$ is \textcolor{darkgreen}{active}, and $m_j = 1$ if $c_j$ is \textcolor{red}{inactive} (for $1 \leq j \leq k$).  Use the notation $F_{m_1 m_2 \cdots m_k}$ to denote a surface $F$ which has the binary string $m_1 m_2 \cdots m_k$ associated to it.

The chain isomorphism will depend on the reordering defined by $\sigma \in S_k$.  We are supposing $\sigma =$ ($i$ $i$+1), and so we define the $\psi_\sigma$ by sending each $D^k$-surface $F$ from the chain complex $\BNC$ to $(+1)$ or $(-1)$ times the corresponding surface in $\BNC'$.  The surface $F$ gets sent to the surface which is the same as $F$ except that the crosscuts $c_i$ and $c_{i+1}$ have their subscript labels interchanged.  Whether the weight $(+1)$ or $(-1)$ is used depends on whether the crosscuts $c_i$ and $c_{i+1}$ are \textcolor{darkgreen}{active} or \textcolor{red}{inactive}.  The weight $(-1)$ is used if both crosscuts are \textcolor{red}{inactive}, and weight $(+1)$ is used otherwise.  Using the notation described in the previous paragraph, we have that

$$
\psi_\sigma(F_{m_1 \cdots m_i m_{i+1} \cdots m_k}) = \left\{ \begin{tabular}{ll} $-F_{m_1 \cdots m_{i+1} m_i \cdots m_k}$ \ & \ if \ $m_i = m_{i+1} = 1$\\ $F_{m_1 \cdots m_{i+1} m_i \cdots m_k}$ \ & \ otherwise \end{tabular} \right.
$$

It is clear that $\psi_\sigma$ is an isomorphism, but we must check that it is a chain map.  To show that $\psi_\sigma$ commutes with the differential $d$, it suffices to show that $\psi_\sigma \circ (-1)^{\alpha(c_j)}d_{c_j} = (-1)^{\alpha(c_{\sigma(j)})}d_{c_{\sigma(j)}} \circ \psi_\sigma$ for $1 \leq j \leq k$, where $\alpha(c)$ is the number of \textcolor{red}{inactive} crosscuts that come before $c$ in the ordering of crosscuts on $F$ (as described in definition \ref{4.3}).  This problem may be seen as showing that a diagram commutes.  The diagram in question is given below.  For ease of notation in constructing the diagram below, we use $1<i<i+1<j<k$.  We will \underline{not} assume that this true in general.

$$
F_{m_1 \cdots m_i m_{i+1} \cdots  0   \cdots m_k} \ \ \ \ \ \ \ \ \ \ \ \ \ \ \ \ \ \ \ \ \   F_{m_1 \cdots m_{i+1} m_i \cdots 0 \cdots m_k}
\put(-146,1){\vector(1,0){53}}
\put(-197,-18){\vector(0,-1){40}}
\put(-43,-18){\vector(0,-1){40}}
\put(-128,9){ \footnotesize  $\psi_\sigma$ }
\put(-251,-38){ \footnotesize $(-1)^{\alpha(c_j)}d_{c_j}$ }
\put(-41,-38){ \footnotesize $(-1)^{\alpha(c_{\sigma(j)})}d_{c_{\sigma(j)}}$ }
$$

\vspace{-6 pt}
$$
F_{m_1 \cdots m_i m_{i+1} \cdots  1   \cdots m_k}  \ \ \ \ \ \ \ \ \ \ \ \ \ \ \ \ \ \ \ \ \  F_{m_1 \cdots m_{i+1} m_i \cdots  1   \cdots m_k}
\put(-146,1){\vector(1,0){53}}
\put(-128,9){ \footnotesize  $\psi_\sigma$ }
$$

\bigskip

From how $\psi$ and $d_{c_j}$ are defined, we know that $\psi_\sigma \circ (-1)^{\alpha(c_j)}d_{c_j} = \pm[(-1)^{\alpha(c_{\sigma(j)})}d_{c_{\sigma(j)}} \circ \psi_\sigma]$.  To show that we always have $\psi_\sigma \circ (-1)^{\alpha(c_j)}d_{c_j} = +[(-1)^{\alpha(c_{\sigma(j)})}d_{c_{\sigma(j)}} \circ \psi_\sigma]$, it suffices to show that

\begin{center}
sign$[\psi_\sigma(F_{m_1 \cdots m_{j-1} \ 0 \ m_{j+1} \cdots m_k})] \neq$ sign$[\psi_\sigma(F_{m_1 \cdots m_{j-1} \ 1 \ m_{j+1} \cdots m_k})]$\\
$ \Longleftrightarrow $\\
sign$[(-1)^{\alpha(c_j)}d_{c_j}] \neq$ sign$[ (-1)^{\alpha(c_{\sigma(j)})}d_{c_{\sigma(j)}}]$.
\end{center}

\begin{itemize}
\item[($\Rightarrow$)] In this case, exactly one of $\psi_\sigma(F_{m_1 \cdots m_{j-1} \ 0 \ m_{j+1} \cdots m_k})$ and $\psi_\sigma(F_{m_1 \cdots m_{j-1} \ 1 \ m_{j+1} \cdots m_k})$ is negative.  Since $\psi_\sigma$ is only negative when $i=i+1=1$, it must be that one of $m_i$ and $m_{i+1}$ is $m_j$ and the other equals $1$.  We then have that $\alpha(c_i) = \alpha(c_{i+1}) \pm 1$.  Hence, $[(-1)^{\alpha(c_j)}d_{c_j}] = - [ (-1)^{\alpha(c_{\sigma(j)})}d_{c_{\sigma(j)}}]$.
\item[($\Leftarrow$)] If sign$[(-1)^{\alpha(c_j)}d_{c_j}] \neq$ sign$[ (-1)^{\alpha(c_{\sigma(j)})}d_{c_{\sigma(j)}}]$, then the number of $1$'s which come before $m_j$ must be different than the number of $1$'s which come before $m_{\sigma(j)}$.  This could not be the case if either $j<i$ or $i+1<j$, so it must be that $j = i$ or $j = i+1$.  Also, if one of $m_i$ and $m_{i+1}$ is $m_j$, and the other equals $0$, then $[(-1)^{\alpha(c_j)}d_{c_j}] = [ (-1)^{\alpha(c_{\sigma(j)})}d_{c_{\sigma(j)}}]$, a contradiction.  Therefore it must be that one of $m_i$ and $m_{i+1}$ is $m_j$, and the other equals $1$.  We then have that $\psi_\sigma(F_{m_1 \cdots m_{j-1} \ 0 \ m_{j+1} \cdots m_k})$ is the positive identity map, and that $\psi_\sigma(F_{m_1 \cdots m_{j-1} \ 1 \ m_{j+1} \cdots m_k})$ is multiplication by $(-1)$. \qed
\end{itemize} 

Recall that in Definition \ref{complex_defn}, we defined the Bar-Natan chain modules $\BNC_{i,j,k,b}(L)$ of the Bar-Natan chain complex $\BNC$ to be the free module of equivalence classes of isotopy classes of $D^k$-surfaces in $S^3$ with $I=i, J=j, K=k,B=b,$ and $\partial(F) = L$.  In other words, the complex $\BNC$ is a sequence of free modules $\BNC_{i,j,k,b}(L)$ with bases that consist of (equivalence classes of isotopy classes of) $D^k$-surfaces.

\begin{defn}
Let $\hspace{1pt} \BNC$ be a Bar-Natan chain complex whose sequence of Bar-Natan modules $\BNC_{i,j,k,b}$ have bases $\{ F^{i,j,k,b}_a \}_{a \in A_{i,j,k,b} }$ and let $B\mathcal{D} \subseteq \BNC$ be a (Bar-Natan) subcomplex of $\BNC$ with bases $\{ F^{i,j,k,b}_a \}_{a \in A'_{i,j,k,b} }$ (where all $A'_{i,j,k,b}$ and $A_{i,j,k,b}$ are indexing sets with $A'_{i,j,k,b} \subseteq A_{i,j,k,b}$ for all $i,j,k,b$).  Use $B\mathcal{D}'$ to denote the complementary (Bar-Natan) subcomplex of $\BNC$ generated by the basis $\{ F^{i,j,k,b}_a \}_{a \in A_{i,j,k,b} - A'_{i,j,k,b} }$.  We say that $B\mathcal{D}$ is an \textit{isolated} subcomplex of $\BNC$ if there are no non-zero local maps $d_c$ between basis elements in $B\mathcal{D}$ and basis elements in $B\mathcal{D}'$.  If $B\mathcal{D}'$ is non-empty, then it follows that $B\mathcal{D}'$ is an isolated subcomplex of $\BNC$ as well.
\end{defn}

\begin{defn}
Given an isolated subcomplex $B\mathcal{D} \subseteq \BNC$, we say that $B\mathcal{D}$ is \textit{i-reducible} if there exists an isolated (nonzero) subcomplex $B\mathcal{D'} \subset \BNC$ with $B\mathcal{D'} \subsetneq B\mathcal{D}$.  Otherwise, we say that $B\mathcal{D}$ is \textit{i-irreducible}.
\end{defn}

The following proposition and corollary are both easy to show -- their proofs are left to the reader.

\begin{prop}
If $B\mathcal{D}_1, B\mathcal{D}_2 \subseteq \BNC$ are isolated subcomplexes, then either $(B\mathcal{D}_1$ and $B\mathcal{D}_2)$ are disjoint or $(B\mathcal{D}_1 \cap B\mathcal{D}_2)$ is an isolated subcomplex of $\BNC$.
\end{prop}

\begin{cor}
All i-irreducible (isolated) subcomplexes of a Bar-Natan chain complex are pairwise disjoint.
\end{cor}

Now we define an equivalence relation on the i-irreducible subcomplexes of a Bar-Natan chain complex $\BNC$ via group action.  Given two i-irreducible subcomplexes $B\mathcal{D}, B\mathcal{D'} \subseteq \BNC$, we say $B\mathcal{D} \sim B\mathcal{D'}$ if and only if there exists a chain isomorphism $\psi_\sigma: \BNC \rightarrow \BNC '$ induced by an action of the symmetric group $S_k$ on the underlying ordering of the crosscuts of the $D^k$-surfaces in $\BNC$ (as in the proof of Theorem \ref{ccthm1}) which maps $B\mathcal{D}$ onto $B\mathcal{D}'$.  

\begin{defn}
Consider the equivalence classes given by the orbits of the group action described above.  Each equivalence class consists of some number (at most $k!$) of i-irreducible subcomplexes of $\BNC$ which only differ by permutations on the ordering of the crosscuts.  Choose one representative (an i-irreducible subcomplex of $\BNC$) for each equivalence class; such a choice corresponds to picking a fixed ordering of the crosscuts.  Taking the union of the bases for these representative subcomplexes creates a basis for new chain complex.  Let this new chain complex be the \textit{diagramless complex} of the link $L$, denoted by $D\mathcal{C}(L)$.  Refer to the resulting homology as the \textit{diagramless homology} of $L$, and denote it by $D\mathcal{H}(L)$.
\end{defn}

\begin{example}
Below we give examples of three i-irreducible (isolated) subcomplexes.  The first two are in the same equivalence class because they only differ by a reordering of the crosscuts.  The third is in a different equivalence class because the crosscut orientations (directions) do not all match up correctly.  In general, equivalence of i-irreducible subcomplexes may be obstructed by differing crosscut orientation, by differing number/placement of dots or by differing number/type of surfaces involved.\\
\begin{center}
An i-irreducible subcomplex (example 1 of 3):\\
	\ \includegraphics[scale=1.5]{u0.pdf}
		\color{black}
		\put(5,27){\vector(1,0){52}}
		\put(25,31){$d_{c_1}$}
		\thicklines
		\color{red}
		\put(-40.2,17){\vector(1,0){34.5}}
		\put(-50.6,13){$c_2$}
		\color{darkgreen}
		\put(-40.2,37){\vector(1,0){34.5}}
		\put(-51,35){$c_1$}
		\color{black}
	\ \ \ \ \ \ \ \ \ \ \ \ \ \ \ \ \ \ \ 
	\includegraphics[scale=1.5]{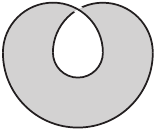} 
		\thicklines
		\color{red}
		\put(-26.2,39){\vector(1,1){12.3}}
		\put(-35.7,34.5){$c_1$}
		\color{red}
		\put(-71,9){$c_2$}
		\put(-63,14){\vector(1,0){52}}
		\color{black}
\end{center}

\begin{center}
An i-irreducible subcomplex (example 2 of 3):\\
	\ \includegraphics[scale=1.5]{u0.pdf}
		\color{black}
		\put(5,27){\vector(1,0){52}}
		\put(25,31){$d_{c_2}$}
		\thicklines
		\color{red}
		\put(-40.2,17){\vector(1,0){34.5}}
		\put(-50.6,13){$c_1$}
		\color{darkgreen}
		\put(-40.2,37){\vector(1,0){34.5}}
		\put(-51,35){$c_2$}
		\color{black}
	\ \ \ \ \ \ \ \ \ \ \ \ \ \ \ \ \ \ \ 
	\includegraphics[scale=1.5]{u1.pdf} 
		\thicklines
		\color{red}
		\put(-26.2,39){\vector(1,1){12.3}}
		\put(-35.7,34.5){$c_2$}
		\color{red}
		\put(-71,9){$c_1$}
		\put(-63,14){\vector(1,0){52}}
		\color{black}
\end{center}

\begin{center}
An i-irreducible subcomplex (example 3 of 3):\\
	\ \includegraphics[scale=1.5]{u0.pdf}
		\color{black}
		\put(5,27){\vector(1,0){52}}
		\put(25,31){$d_{c_2}$}
		\thicklines
		\color{red}
		\put(-5.7,17){\vector(-1,0){34.5}}
		\put(-50.6,13){$c_1$}
		\color{darkgreen}
		\put(-40.2,37){\vector(1,0){34.5}}
		\put(-51,35){$c_2$}
		\color{black}
	\ \ \ \ \ \ \ \ \ \ \ \ \ \ \ \ \ \ \ 
	\includegraphics[scale=1.5]{u1.pdf} 
		\thicklines
		\color{red}
		\put(-26.2,39){\vector(1,1){12.3}}
		\put(-35.7,34.5){$c_2$}
		\color{red}
		\put(-71,9){$c_1$}
		\put(-11,14){\vector(-1,0){52}}
		\color{black}
\end{center}

\end{example}

\section{An Injection \texorpdfstring{$\mathcal{K}_k(D) \hookrightarrow \mathcal{C}_{i,j,k,b}(L)$}{}}
\label{6}

In \cite{Kh1} Mikhail Khovanov defined a homology for a given diagram of a link.  In this section it is shown that for a given link diagram $D$ of a link $L$, there is an injection from the chain complex for Khovanov homology of $D$ into the chain complex for the diagramless homology of $L$.   We obtain an injective chain map $\iota$ by finding chain modules in the diagramless homology that correspond to the chain modules of Khovanov homology.  To do this, we must first develop language to talk about the Khovanov chain modules.  

Recall that in his state sum for the Jones polynomial Kauffman \cite{Kau} defines a \textit{state} of a link diagram to be a collection of markers (one for each crossing) that specify a pair of opposite angles.  A marker at a crossing defines a `smoothing' of that crossing, depending on the the type of marker present.  

\begin{center}
\includegraphics[scale = 1]{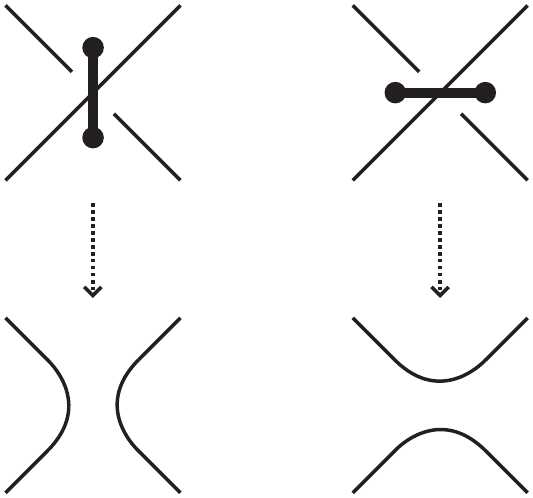}
\put(-162.6,150){ {\small positive marker} }
\put(-62.6,150){ {\small negative marker} }
\put(-168,-10){ {\small positive smoothing} }
\put(-68,-10){ {\small negative smoothing} }
\end{center}

A state $s_D$ of a link diagram $D$ defines a smoothing at each crossing.  Hence to each state $s_D$ there is a planar collection of disjoint (possible embedded) circles obtained by smoothing each crossing of the link diagram -- this is referred to as the `complete smoothing'.

We borrow notation from Viro's paper \cite{Vi} and use \textit{enhanced Kauffman state} to refer to a state $s_D$ which is enhanced by an assignment of a plus or minus to each circle in the complete smoothing of $s_D$.  A capital $S_D$ is used to refer to an enhanced Kauffman state.  

We define the chain module $\mathcal{K}_k$ for Khovanov homology to be the free module of isotopy classes of enhanced Kauffman states with exactly $k$ positive markers.  The differential for the Khovanov chain complex can be defined as a certain combination of local maps which consist of replacing a positive marker by a negative marker.  This differential and these local maps will be defined later in this section.

Each Khovanov chain module is generated by enhanced Kauffman states with a specified number of positive markers.  Hence, to define the chain map $\iota: \mathcal{K}_k(D) \hookrightarrow \mathcal{C}_{i,j,k,b}(L)$ on chain modules of Khovanov homology, it suffices to show where $\iota$ sends enhanced Kauffman states.  Given an enhanced Kauffman state $S_D$, we will define a $D^k$-surface $F_{S_D}$ and set $\iota(S_D):=F_{S_D}$.   We will call $F_{S_D}$ the \textit{state surface} corresponding to $S_D$.\footnote{In \cite{Oz}, Ozawa uses the name $\sigma$-\textit{state surface} for something similar to our state surfaces.  The `$\sigma$-' prefix includes the information of the Kauffman state from which the surface is derived.  Our state surfaces will be built in a similar fashion, but from \textit{enhanced} Kauffman states.}

\begin{defn}
\label{statesurface}
Let $D$ be a link diagram with a given ordering of its crossings and let $S_D$ an enhanced Kauffman state of $D$.  The \textit{state surface} corresponding to $S_D$ is denoted by $F_{S_D}$ and is built from $S_D$ as follows:

Consider the complete smoothing of the enhanced Kauffman state $S_D$.  The complete smoothing of $S_D$ consists of a finite number of disjoint (possibly embedded) circles each marked with a plus or a minus.  Assume there exists an oriented 2-sphere $\Sigma$ embedded in $S^3$ and there exist closed 3-balls $B^3_+$ and $B^3_-$ embedded in $S^3$ such that the complete smoothing of $S_D$ is in $\Sigma$, and $\Sigma = B^3_+ \cap B^3_-$.

\begin{center}
	\includegraphics[trim = 0mm 46mm 0mm 0mm, clip,scale=1]{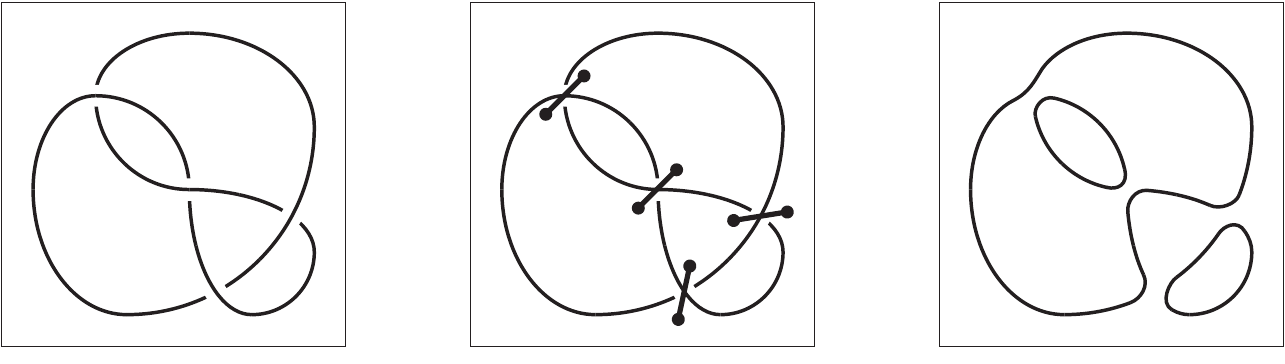} 
\put(-361,-6){A link diagram $D$}
\put(-235,-6){An enhanced state $S_D$}
\put(-198,60){$-$}
\put(-161,22){$+$}
\put(-214,23){$+$}
\put(-105,-6){The complete smoothing}
\put(-62.5,60){$-$}
\put(-24.5,22.5){$+$}
\put(-77.4,22){$+$}
\end{center}

For each circle of the complete smoothing, take a copy of the disk it bounds in $\Sigma$ and push the interior of the disk into $B^3_+$.  The disk interiors should be pushed into $B^3_+$ so that they do not interest one another.  For each disk, if the disk is bounded by a circle marked with a plus, place a dot on the interior of the disk.  Disks bounded by circles marked by a minus receive no dots.  At this point, the plus/minus information of the circles can be forgotten.

\begin{center}
	\includegraphics[trim = 0mm 0mm 0mm 47mm, clip,scale=1]{statesurf_1.pdf} 
\put(-381,-11){Creating disks from circles}
\put(-294.5,17.5){\textcolor{blue}{$\bullet$}}
\put(-348,22){\textcolor{blue}{$\bullet$}}
\put(-235,-11){The cross-dual surface}
\put(-159.5,17.5){\textcolor{blue}{$\bullet$}}
\put(-213,22){\textcolor{blue}{$\bullet$}}
\put(-88,-11){The state surface}
\put(-24,18){\textcolor{blue}{$\bullet$}}
\put(-77.5,22){\textcolor{blue}{$\bullet$}}
\end{center}

To obtain the cross-dual surface, insert a locally oriented rectangular strip of surface in $\Sigma$ that connects two disks at each location of a former crossing.  The locally oriented pieces should be oriented so that they appear to be positively oriented when viewed from the $B^3_-$ side of $\Sigma$.

To obtain the the \textit{state surface}, $F_{S_D}$, replace each locally oriented rectangular strip of surface by a (non-oriented) piece of surface with a $\frac{1}{2}$-twist that has a crosscut running across it.  If the site of the former crossing (of the Kauffman state) had a positive marker, a right-handed $\frac{1}{2}$-twist with an \textcolor{darkgreen}{active} crosscut should be placed.  If the site of the former crossing (of the Kauffman state) had a negative marker, a left-handed $\frac{1}{2}$-twist with an \textcolor{red}{inactive} crosscut should be placed.  In either case, the crosscut should be oriented to point away from the $B^3_+$ region towards the $B^3_-$ region\footnote{As stated above, the crosscuts placed on the $\frac{1}{2}$-twisted pieces of surface should be oriented to point away from the $B^3_+$ region towards the $B^3_-$ region.  Since the surface is depicted as being viewed from the $B^3_-$ side of $\Sigma$, all crosscuts should appear to point toward us.  Also, this is equivalent to orienting the crosscuts in such a way that replacing crosscut neighborhoods to obtain the cross-dual (as in Definition \ref{2.6}) would result in the a cross-dual surface with the same local orientations as was obtained in the previous step.}, and the boundary of the $\frac{1}{2}$-twisted piece of surface inserted should agree with the original link diagram.  By construction, the number of crossings of $D$ is equal to number of crosscuts of $F_{S_D}$.  Use the ordering of the crossings of $D$ to induce the ordering of the crosscuts of $F_{S_D}$. 
\end{defn}

\begin{remark}
In Section \ref{dhsection} the diagramless chain complex was obtained from the Bar-Natan chain complex by identifying all subcomplexes that differed only by a reordering of the crosscuts.  With the diagramless complex, instead of considering all $k!$ different ways to orient $k$ crosscuts, we can pick an arbitrary fixed ordering of the crosscuts that is consistent between all related $D^k$-surfaces.  A similar approach can be taken when constructing the Khovanov chain complex of a link diagram; an arbitrary fixed ordering of the crossings is chosen.  

In this section, an arbitrary ordering of the crossings of the link diagram is used.  This ordering of the diagram's crossings induces an ordering of the crosscuts of the resulting state surfaces, as we just saw in Definition \ref{statesurface}.
\end{remark}

Given an enhanced Kauffman state $S_D$, we wish to define $\iota(S_D):=F_{S_D}$.  To do this, it must be shown that the state surface $F_{S_D}$ is a $D^k$-surface, the proof of which is somewhat technical.

\begin{prop}
\label{state_surf_prop}
Given an enhanced Kauffman state $S_D$, the resulting state surface $F_{S_D}$ is a $D^k$-surface.
\end{prop}

\emph{Proof.}  In definition \ref{statesurface}, the state surface (as well at its cross-dual and skeleton) are built to satisfy the conditions required to be a $D^k$-surface.  The disk interiors involved in the construction, which are the facets of the cross-dual surface, live in $B^3_+$.  The rest of the cross-dual surface (the skeleton of the cross-dual) is in the embedded oriented 2-sphere $\Sigma$.  The locally oriented pieces of surface are all oriented to agree with the orientation of $\Sigma$.  The only $D^k$-surface condition that is not immediately seen to be satisfied from the construction of a state surface is the orientability of the cross-dual.  

We now show that the cross-dual of the state surface is orientable.  Consider the underlying `complete smoothing' for the state surface.  The complete smoothing is a set of (possibly embedded) circles in $\Sigma$.  Pick a point in $p \in \Sigma$ that lies outside of the original link diagram.  The label each circle bounding that region with a `1'.  Remove these 1-circles from this region, then mark each circle bounding this new region with a `2'.  Again, remove the 2-circles and mark the next set of bounding circles (if any) with a `3'.  Continue this process until all circles are marked with a positive integer.  

These circles separate $\Sigma$ into regions.  To the region containing the point $p$, give a (temporary) local orientation that agrees with the orientation of $\Sigma$.  To all regions between 1-circles and 2-circles give (temporary) local orientations that disagree with that of $\Sigma$.  Continue in this fashion, giving successive regions alternating orientations, until all such regions have a (temporary) local orientation.  

Since these (temporary) local orientations of the regions alternate, they induce a (temporary) orientation on the circles.  In turn, we can use the orientations of the circles to induce orientations on the facets (disks) that they bound in the cross-dual surface.  Give the facets of the cross-dual a (permanent) orientation that agrees with the orientations of their boundary circles.  An example is given below.
\begin{center}
	\includegraphics[trim = 95mm 47.5mm 0mm 0mm, clip,scale=1]{statesurf_1.pdf} 
\put(-50,61){{\footnotesize $\leftarrow{2}$-circle}}
\put(-8.5,24){{\footnotesize $\longleftarrow{1}$-circle}}
\put(-11.5,74){{\footnotesize $\longleftarrow{1}$-circle}}
\ \ \ \ \ \ \ \ \ \ \ \ 
	\includegraphics[trim = 95mm 47.5mm 0mm 0mm, clip,scale=1]{statesurf_1.pdf} 
\put(-50,61){{\footnotesize $\leftarrow{2}$-circle}}
\put(-8.5,24){{\footnotesize $\longleftarrow{1}$-circle}}
\put(-11.5,74){{\footnotesize $\longleftarrow{1}$-circle}}
\put(-66,57){ \textbf{$\circlearrowleft$} }
\put(-28,18){ \textbf{$\circlearrowright$} }
\put(-47,78){ \textbf{$\circlearrowright$} }
\put(-93,84){ \textbf{$\circlearrowleft$} }
\ \ \ \ \ \ \ \ \ \ \ \ 
	\includegraphics[trim = 47mm 0mm 47mm 47.5mm, clip,scale=1]{statesurf_1.pdf} 
\put(-68,57){ \textbf{$\circlearrowleft$} }
\put(-30,18){ \textbf{$\circlearrowright$} }
\put(-45,73){ \textbf{$\circlearrowright$} }
\put(-45,13.5){ \textcolor{lightblue}{$\bullet$} }
\put(-80.5,70){ \textcolor{lightblue}{$\bullet$} }
\put(-53,43){ \textcolor{lightblue}{$\bullet$} }
\put(-24,35.7){ \textcolor{lightblue}{$\bullet$} }
\end{center}
To give a (global) orientation to the cross-dual surface, we let the pieces of surface that connect the oriented disks receive an orientation consistent with those disks.  However, we must check that such an orientation is possible. 

First, some terminology.  Let $\sigma$ be one of the circles from the complete smoothing.  Use the terms \textit{outside} and \textit{inside} of $\sigma$ to mean the regions of $(\Sigma - \sigma)$ that contain the point $p$ and do not contain the point $p$, respectively.   Since consecutive embedded circles have consecutive integer labels, the pieces of surface connecting the circles in $\Sigma$ either connect circles with the same label or connect circles with labels that differ by plus or minus one.  

If a piece of surface connects circles with the same label, the piece of surface is either connecting a circle to itself or is connecting neighboring circles (circles that do not appear to be embedded in one another with respect to the outside region containing the point $p$).  In the former case, the piece of surface is either on the inside or outside of the circle, meaning it can be oriented to disagree with the circle's local orientation at both spots.  In the latter case, the piece of surface is on the outside of both circles, connecting them.  Since those two circles have the same label, the piece of surface may be oriented to disagree with both of them.

If a piece of surface connects two circles with different labels, it must be connecting a circle embedded inside of another (with respect to the outside region containing the point $p$).  Hence, the piece of surface connects to the outside of one circle and to the inside of the other.  This implies that the piece of surface may be oriented to disagree with the local orientation of both circles (since the circles have opposite types of local orientations based on their different labels).  

Therefore, the cross-dual surface of the state surface is orientable.  Since all conditions were met, the state surface is a $D^k$-surface. \qed
\bigskip

Next we describe how $\iota$ relates the differential of the Khovanov chain complex to differential of the diagramless chain complex.  Recall that in Definition \ref{4.3} the differential $d$ for the diagramless chain complex is defined to be the weighted sum of local operators $d_c$, where $d_c$ acts on a neighborhood of the crosscut $c$.  Similarly, the differential $\partial$ for the Khovanov complex can be defined as a weighted sum of certain local operators.  We will use the notation $\partial_i$ to refer to the local operator that acts on the $i^{\textrm{th}}$ crossing of the link diagram D.  

For the Khovanov chain complex, the local operators $\partial_i$ act on an enhanced Kauffman state $S_D$ by replacing a positive marker with a negative marker, where the crossing involved was labeled as the $i^{\textrm{th}}$ crossing of $D$.  This is equivalent to changing a positively smoothed crossing into a negatively smoothed one.  This process either results in the merging of two distinct circles or results in the splitting of one circle into two.  The fact that circles are marked with a plus or minus sign in (smoothed) enhanced Kauffman states complicates the situation.  

Due to Khovanov's differential $\partial$ being bidegree (1,0), there are certain restrictions on the $\pm$ markings of the resulting circle(s) after merging or splitting.\footnote{The Khovanov chain complex can be seen to be bigraded with a homological grading and a polynomial grading.  In this paper, the polynomial grading of the Khovanov complex is ignored.  However, the chain complex for the diagramless homology is bestowed with a polynomial grading which is designed to have a natural correspondence with the Khovanov polynomial grading.}  In Viro's paper \cite{Vi} he shows that these restrictions have the following implications for the local operator $\partial_i$:

\begin{enumerate}
\item If applying $\partial_i$ corresponds to the splitting of one positively marked circle, then both of the resulting circles are positively marked.
\begin{center}
\includegraphics[trim = 78mm 0mm 0mm 0mm,clip,scale=1]{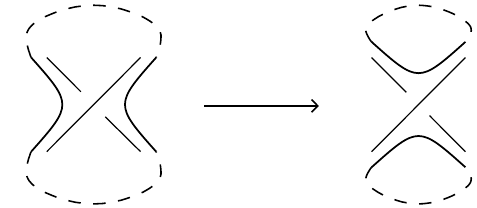}
\put(-115,49.5){+}
\put(-20.5,49.5){+}
\put(-20.5,5.5){+}
\put(-66,33){$\partial_i$}
\end{center}
\item If applying $\partial_i$ corresponds to the splitting of one negatively marked circle, then the result is a sum of two different enhanced Kauffman states, namely the ones obtained from the two different ways to positively mark one circle and negatively mark the other.
\begin{center}
\includegraphics[trim = 78mm 0mm 0mm 0mm,clip,scale=1]{vpic.pdf}
\put(-115,49.5){$-$}
\put(-20.5,49.5){+}
\put(-20.5,5.5){$-$}
\put(-66,33){$\partial_i$}
\ \ \ \ \ 
\includegraphics[trim = 113mm 0mm 0mm 0mm,clip,scale=1]{vpic.pdf}
\put(-50,27){{\Large \textbf{+}}}
\put(-20.5,49.5){$-$}
\put(-20.5,5.5){+}
\end{center}
\item If applying $\partial_i$ corresponds to the merging of two circles that are both negatively marked, then the resulting circle is also negatively marked.
\begin{center}
\includegraphics[trim = 0mm 0mm 64mm 0mm,clip,scale=1]{vpic.pdf}
\put(-173,27.5){$-$}
\put(-129.5,27.5){$-$}
\put(-18.5,27.5){$-$}
\put(-95.5,33){$\partial_i$}
\end{center}
\item If applying $\partial_i$ corresponds to the merging of two circles, where one is negatively marked and the other is positively marked, then the resulting circle is positively marked.
\begin{center}
\includegraphics[trim = 0mm 0mm 64mm 0mm,clip,scale=1]{vpic.pdf}
\put(-173,27.5){$-$}
\put(-129.5,27.5){+}
\put(-18.5,27.5){+}
\put(-95.5,33){$\partial_i$}
\end{center}
\item If applying $\partial_i$ corresponds to the merging of two circles that are both positively marked, then $\partial_i$ is the zero-map.
\begin{center}
\includegraphics[trim = 0mm 0mm 88mm 0mm,clip,scale=1]{vpic.pdf}
\put(-105,27.5){+}
\put(-61.5,27.5){+}
\put(-27.5,33){$\partial_i$}
\put(1,23){{\huge 0}}
\put(15,26){{\large (zero)}}
\ \ \ \ \ \ \ \ \ \ \ \ \ \ \ \ \ \ \ \ 
\end{center}
\end{enumerate}

Now we use these local operators to define the Khovanov chain complex differential $\partial$, which takes an enhanced Kauffman state of a link diagram $D$ with $n$ positively smoothed crossings to an enhanced Kauffman state of $D$ with $n-1$ positively smoothed crossings.

$$
\partial := \sum_{\begin{array}{c} \textrm{enhanced Kauffman states} \\ S_D \textrm{ of } D \end{array}} \sum_{\begin{array}{c} \textrm{positively smoothed} \\ \textrm{crossings of } D \end{array}} (-1)^{\sigma(i)} \partial_i
$$
The crossings of $D$ are ordered $1,...,k$, and $\sigma(i)$ is the number positively smoothed crossings of $S_D$ that are come before $i$ in the ordering $1,...,k$.

Now we are ready to describe how $\iota$ relates the differential of the Khovanov chain complex to differential of the chain complex for the diagramless homology.  

Let $L$ be a link, and $D$ be a link diagram of $L$ with crossings $1,...,k$.  Since $\partial$ is defined in terms of the $\partial_i$, it suffices to say where $\iota$ sends the $\partial_i$.  Define $\iota(\partial_i)$ to be the local map $d_{c_i}$, where $d_{c_i}$ acts on the $i^{\textrm{th}}$ crosscut of a $D^k$-surface with boundary $L$, as described in Definition \ref{4.1}.
$$
\iota(\partial_i) = d_{c_i}
$$
In order for $\iota$ to be a chain map from the Khovanov chain complex to the chain complex for diagramless homology, it must satisfy the relation $d_{c_i} \hspace{-2 pt}\circ \iota = \iota \circ \partial_i$.  Since the chain (modules/groups) for the Khovanov complex are generated by enhanced Kauffman states, it suffices to show that this relation holds for an arbitrary enhanced Kauffman state.

Let $S_D$ be an enhanced Kauffman state of a link diagram $D$.  We must show that $d_{c_i}\hspace{-2 pt}\left( \iota (S_D) \right) = \iota \left( \partial_i (S_D) \right)$.
By definition of $\iota$, $d_{c_i}\hspace{-2 pt}\left( \iota (S_D) \right) = \iota \left( \partial_i (S_D) \right)$ now becomes $d_{c_i}\hspace{-2 pt}\left( F_{S_D} \right) = F_{\partial_i (S_D)}$.  Hence, we must show that applying $d_{c_i}$ to the state surface for $S_D$ gives the state surface for the enhanced Kauffman state $F_{\partial_i (S_D)}$.  The five different situations of the local operator $\partial_i$ will be considered separately.
\begin{enumerate}
\item This case involves applying $\partial_i$ when it corresponds to the splitting of one positively marked circle into two positively marked circles.  Applying $\iota$ to the result gives the state surface $F_{\partial_i (S_D)}$, which has 1 dot on each of the two neighboring facets of the crosscut $c_i$.  On the other hand, applying $\iota$ first gives the state surface $F_{S_D}$, which has only one facet.  This facet has 1 dot.  Applying $d_{c_i}$ reveals the presence of a compressing disk.  Using the relation (NC) to compress and then using (D2) yields the same surface as before, as seen in the figure below.  
\begin{center}
\begin{tabular}{c}
{
\includegraphics[trim = 78mm 0mm 0mm 0mm,clip,scale=1]{vpic.pdf}
\put(-115,49.5){+}
\put(-20.5,49.5){+}
\put(-20.5,5.5){+}
\put(-66,33){$\partial_i$}
\put(-145,-5){\vector(-1,-1){28}}
\put(-165,-18){$\iota$}
\put(7,-6){\vector(1,-1){28}}
\put(24,-19){$\iota$}
\put(-72,-25){{\Huge $\circlearrowright$}}
} \\
 \ \\
\reflectbox{ 
\includegraphics[trim = 35.8mm 0mm 0mm 0mm,clip,scale=1,angle=90]{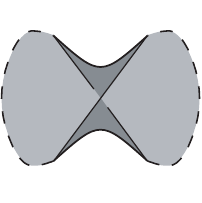}
}
\put(-39,48){\textcolor{blue}{$\bullet$}}
\put(-6,29){\vector(1,0){24}}
\put(-1,34){$d_{c_i}$}
\color{darkgreen}
\put(-24.8,23.7){\vector(0,1){11}}
\color{black}
\ \ \ \ \ \ 
\includegraphics[trim = 20.8mm 0mm 20.8mm 0mm,clip,scale=1]{vsurf2.pdf} 
\put(-29,48){\textcolor{blue}{$\bullet$}}
\color{red}
\put(-26,33){\vector(1,0){10}}
\color{black}
\ \ \ \ \ 
\reflectbox{ \includegraphics[trim = 0mm 0mm 35.8mm 0mm,clip,scale=1,angle=90]{vsurf2.pdf} }
\put(-34,42){\textcolor{blue}{$\bullet$}}
\put(-39,50){\textcolor{blue}{$\bullet$}}
\put(-77,26){{\Large \textbf{$\stackrel{(\textrm{NC})}{=}$}}}
\color{red}
\put(-37.3,36){\vector(1,0){10}}
\color{black}
\hspace{-30 pt}
\reflectbox{ \includegraphics[trim = 0mm 0mm 35.8mm 0mm,clip,scale=1,angle=90]{vsurf2.pdf} }
\put(-33,10){\textcolor{blue}{$\bullet$}}
\put(-38,50){\textcolor{blue}{$\bullet$}}
\put(-58,26){{\Large \textbf{+}}}
\color{red}
\put(-37.3,36){\vector(1,0){10}}
\color{black}
\reflectbox{ \includegraphics[trim = 0mm 0mm 35.8mm 0mm,clip,scale=1,angle=90]{vsurf2.pdf} }
\put(-33,10){\textcolor{blue}{$\bullet$}}
\put(-38,50){\textcolor{blue}{$\bullet$}}
\put(-75,26){{\Large \textbf{$\stackrel{(\textrm{D2})}{=}$}}}
\color{red}
\put(-37.3,36){\vector(1,0){10}}
\color{black}
\end{tabular}
\end{center}
\item This case is very similar to case 1.  The difference here is that applying $\partial_i$ corresponds to the splitting of one \textit{negatively} marked circle, resulting in a sum of two different enhanced Kauffman states, namely the ones obtained from the two different ways to positively mark one circle and negatively mark the other.  The other difference is the absence of dots on certain facets.  However, just the right number of dots are absent in just the right places.  We get the same result if we apply $\partial_i$ then $\iota$ or if we apply $\iota$ then $d_{c_i}$.
\begin{center}
\begin{tabular}{c}
{
\includegraphics[trim = 78mm 0mm 0mm 0mm,clip,scale=1]{vpic.pdf}
\put(-115,49.5){$-$}
\put(-20.5,49.5){+}
\put(-20.5,5.5){$-$}
\put(-66,33){$\partial_i$}
\ \ \ \ \ 
\includegraphics[trim = 113mm 0mm 0mm 0mm,clip,scale=1]{vpic.pdf}
\put(-50,27){{\Large \textbf{+}}}
\put(-20.5,49.5){$-$}
\put(-20.5,5.5){+}
\put(-174,-10){\vector(-1,-3){8}}
\put(-184,-21){$\iota$}
\put(-37,-10){\vector(1,-3){8}}
\put(-30.5,-21){$\iota$}
\put(-128,-23){{\huge $\circlearrowright$}}
} \\
 \ \\
\reflectbox{ 
\includegraphics[trim = 35.8mm 0mm 0mm 0mm,clip,scale=1,angle=90]{vsurf2.pdf}
}
\put(-6,29){\vector(1,0){24}}
\put(-1,34){$d_{c_i}$}
\color{darkgreen}
\put(-24.8,23.7){\vector(0,1){11}}
\color{black}
\ \ \ \ \ \ 
\includegraphics[trim = 20.8mm 0mm 20.8mm 0mm,clip,scale=1]{vsurf2.pdf} 
\color{red}
\put(-29.5,33){\vector(1,0){10}}
\color{black}
\ \ \ \ \ 
\reflectbox{ \includegraphics[trim = 0mm 0mm 35.8mm 0mm,clip,scale=1,angle=90]{vsurf2.pdf} }
\put(-34,42){\textcolor{blue}{$\bullet$}}
\put(-77,26){{\Large \textbf{$\stackrel{(\textrm{NC})}{=}$}}}
\color{red}
\put(-37.3,36){\vector(1,0){10}}
\color{black}
\hspace{-30 pt}
\reflectbox{ \includegraphics[trim = 0mm 0mm 35.8mm 0mm,clip,scale=1,angle=90]{vsurf2.pdf} }
\put(-33,10){\textcolor{blue}{$\bullet$}}
\put(-58,26){{\Large \textbf{+}}}
\color{red}
\put(-37.3,36){\vector(1,0){10}}
\color{black}
\end{tabular}
\end{center}
\item Here applying $\partial_i$ corresponds to the merging of two circles that are both negatively marked, and the resulting circle is also negatively marked.  For this case it is easy to verify that $d_{c_i}\hspace{-2 pt}\left( \iota (S_D) \right) = \iota \left( \partial_i (S_D) \right)$ since no relations are involved.
\begin{center}
\begin{tabular}{c}
{
\includegraphics[trim = 0mm 0mm 64mm 0mm,clip,scale=1]{vpic.pdf}
\put(-173,27.5){$-$}
\put(-129.5,27.5){$-$}
\put(-18.5,27.5){$-$}
\put(-95.5,33){$\partial_i$}
\put(-145,0){\vector(0,-1){33}}
\put(-151.5,-18){$\iota$}
\put(-39,0){\vector(0,-1){33}}
\put(-36,-19){$\iota$}
\put(-100,-19){{\huge $\circlearrowright$}}
} \\
\includegraphics[trim = 0mm 0mm 35.8mm 0mm,clip,scale=1]{vsurf2.pdf}
\put(9.5,29){\vector(1,0){24}}
\put(14.5,34){$d_{c_i}$}
\color{darkgreen}
\put(-22.4,24.4){\vector(0,1){10}}
\color{black}
\ \ \ \ \ \ \ \ \ \ \ \ \ 
\includegraphics[trim = 35.8mm 0mm 0mm 0mm,clip,scale=1]{vsurf2.pdf}
\color{red}
\put(-37.5,37){\vector(1,0){10}}
\color{black}
\end{tabular}
\end{center}
\item Here applying $\partial_i$ corresponds to the merging of two circles, where one is negatively marked and the other is positively marked, and the resulting circle is positively marked.  Like case 3, it is easy to verify that $d_{c_i}\hspace{-2 pt}\left( \iota (S_D) \right) = \iota \left( \partial_i (S_D) \right)$ since no relations are involved.  
\begin{center}
\begin{tabular}{c}
{
\includegraphics[trim = 0mm 0mm 64mm 0mm,clip,scale=1]{vpic.pdf}
\put(-173,27.5){$-$}
\put(-129.5,27.5){+}
\put(-18.5,27.5){+}
\put(-95.5,33){$\partial_i$}
\put(-145,0){\vector(0,-1){33}}
\put(-151.5,-18){$\iota$}
\put(-39,0){\vector(0,-1){33}}
\put(-36,-19){$\iota$}
\put(-100,-19){{\huge $\circlearrowright$}}
} \\
\includegraphics[trim = 0mm 0mm 35.8mm 0mm,clip,scale=1]{vsurf2.pdf}
\put(9.5,29){\vector(1,0){24}}
\put(14.5,34){$d_{c_i}$}
\put(-11,28){\textcolor{blue}{$\bullet$}}
\color{darkgreen}
\put(-22.4,24.4){\vector(0,1){10}}
\color{black}
\ \ \ \ \ \ \ \ \ \ \ \ \ 
\includegraphics[trim = 35.8mm 0mm 0mm 0mm,clip,scale=1]{vsurf2.pdf}
\color{red}
\put(-37.5,37){\vector(1,0){10}}
\color{black}
\put(-16,29){\textcolor{blue}{$\bullet$}}
\end{tabular}
\end{center}
\item Here applying $\partial_i$ corresponds to the merging of two circles that are both positively marked, and $\partial_i$ is defined to be the zero-map.  Since $\iota$ is a homomorphism, $\iota(0)=0$.  Hence, we must have that applying $\iota$ followed by $d_{c_i}$ also gives zero.  In this case, $\iota$ is applied to two positively marked circles and gives a twice dotted state surface.  Applying $d_{c_i}$ to this surface yields a single faceted surface with two dots, which equals zero by the (D2) relation.
\begin{center}
\begin{tabular}{c}
{
\hspace{-60pt} \includegraphics[trim = 0mm 0mm 88mm 0mm,clip,scale=1]{vpic.pdf} 
\put(-105,27.5){+}
\put(-61.5,27.5){+}
\put(-27.5,33){$\partial_i$}
\put(1,23){{\huge 0}}
\put(15,26){{\large (zero)}}
\put(-94,-2){\vector(-1,-2){14}}
\put(-106.5,-15){$\iota$}
\put(37,-1){\vector(1,-2){14}}
\put(46,-14){$\iota$}
\put(-39,-14){{\huge $\circlearrowright$}}
} \\
\hspace{-85pt} \includegraphics[trim = 0mm 0mm 35.8mm 0mm,clip,scale=1]{vsurf2.pdf}
\put(9.5,29){\vector(1,0){24}}
\put(14.5,34){$d_{c_i}$}
\put(-51,28){\textcolor{blue}{$\bullet$}}
\put(-13,28){\textcolor{blue}{$\bullet$}}
\color{darkgreen}
\put(-22.4,24.4){\vector(0,1){10}}
\color{black}
\ \ \ \ \ \ \ \ \ \ \ \ \ 
\includegraphics[trim = 35.8mm 0mm 0mm 0mm,clip,scale=1]{vsurf2.pdf}
\color{red}
\put(-37.5,37){\vector(1,0){10}}
\color{black}
\put(-18,28){\textcolor{blue}{$\bullet$}}
\put(-52,28){\textcolor{blue}{$\bullet$}}
\put(31,23){{\huge 0}}
\put(45,26){{\large (zero)}}
\put(4,26){{\Large \textbf{$\stackrel{(\textrm{D2})}{=}$}}}
\end{tabular}
\end{center}
\end{enumerate}

The above five cases prove the following proposition. 

\begin{prop}
\label{inj_prop1}
The map $\iota$, from the Khovanov complex to the diagramless complex, is a chain map.
\end{prop}

We will now show that $\iota$ is an injective chain map.

\begin{prop}
\label{inj_prop2}
The chain map $\iota$ is injective.
\end{prop}

\emph{Proof.}  We will show that the kernel of $\iota$ is trivial.  Since the Khovanov chain complex is generated by enhanced Kauffman states, it suffices to show that if applying $\iota$ to an enhanced Kauffman state equals zero, then this implies the Kauffman state is zero as well.

Let $S_D$ be an enhanced Kauffman state and suppose $\iota(S_D) = 0$.  By definition of $\iota$, $S_D$ is sent to the corresponding state surface $F_{S_D}$.  Now it will be shown that $F_{S_D}$ equals the zero surface without the use of the relations (S0), (S1), (D2) or (NC).  

The process of building a state surface never produces any spheres (dotted or undotted), and so $F_{S_D}$ cannot equal zero via the relations (S0) and (S1).  Since state surfaces are built from disks connected by bands with crosscuts on them, each facet must be incompressible.  That is, a compressing disk on $F_{S_D}$ would have to run across a crosscut because disks are incompressible surfaces.  The relation (NC) is not allowed if the boundary of the compressing disk intersects a crosscut, thus the relation (NC) cannot be used to replace $F_{S_D}$ with surfaces with additional dots.  Since the process of building $F_{S_D}$ places at most one dot on each facet, and since (NC) cannot be used to place additional dots on any facet, the surface $F_{S_D}$ cannot equal zero via the (D2) relation.  Therefore, $\iota(S_D) = F_{S_D} = 0$ without the use of relations.

Since $F_{S_D}$ equals zero directly, this implies that the corresponding (smoothed) enhanced Kauffman state $S_D$ had zero circles.  That is, $S_D = 0$.  \qed

\bigskip
The fact that $\iota$ is an injective chain map implies that the given Khovanov complex can be embedded in the diagramless complex.  However, there is a much stronger statement that can be made.  In the following section, it is proved that the state surfaces from the embedded Khovanov complex (when considering all possible diagrams for the given link) span the set of all $D^k$-surfaces for the diagramless complex.  In other words, the diagramless complex only consists of embedded copies of the Khovanov complex.

\section{Relating the Diagramless Homology to Khovanov Homology}
\label{newsection}

This section begins with a remark and some lemmas that will be used to help prove Proposition \ref{big_prop}, which says that every $D^k$-surface is equal to a linear combination of state surfaces.  This in turn implies that the diagramless complex is comprised entirely of embedded copies of Khovanov complexes for diagrams of the given link.  In the end we are able to prove Theorem \ref{decomp_thm}, which states that the diagramless homology of a link is equal to the direct sum of some number of copies of Khovanov homology for that link. 

\begin{remark}
The relation (NC) can be used (repeatedly) to compress a $D^k$-surface and replace it with the linear combination of `incompressible' $D^k$-surfaces.  It may not be clear what incompressible means in this situation.  When we use the relation (NC) so that a facet $f$ of a $D^k$-surface $F$ is incompressible, this means that $f$ is incompressible in $f \cup (S^3 - F)$, not necessarily that $f$ is incompressible in $S^3$.  If every facet of a $D^k$-surface is incompressible in this way, we say that the $D^k$-surface is incompressible.
\end{remark}

\begin{lemma}
\label{facets_are_disks}
If $F$ is a $D^k$-surface that does not contain any sphere components and has the property that each facet $f$ of $F$ is incompressible in $f \cup (S^3 - F)$, then all of the facets of $F$ are disks.
\end{lemma}

\emph{Proof.}  Let $F$ be a such a $D^k$-surface.  Consider the cross-dual $\fcd$.  By Lemma \ref{facet_lemma}, the facets of $F$ are isotopic to the facets of $\fcd$, and so $\fcd$ does not contain any sphere components and has the property that each facet $f$ of $F$ is incompressible in $f \cup (S^3 - F)$.  Also by Lemma \ref{facet_lemma}, it suffices to prove that the facets of $\fcd$ are all disks.

By the definition of $D^k$-surface there exists an embedded oriented 2-sphere $\Sigma \subset S^3$ and embedded 3-balls $B^3_+, B^3_- \subseteq S^3$ such that 
\begin{itemize}
\item[(1)] $\Sigma = B^3_+ \cap B^3_-$, 
\item[(2)] skel$(\fcd) \subseteq \Sigma$, 
\item[(3)] all of the locally oriented pieces of surface of $\fcd$ agree with the orientation of $\Sigma$, and 
\item[(4)] $\fcd - \textrm{skel}(\fcd) = \{\textrm{the facets of }\fcd\} \subseteq B^3_+$.  
\end{itemize}
Since the boundary of $\fcd$ is a subset of the skeleton of $\fcd$, and since skel$(\fcd) \subseteq \Sigma$, we know that $\partial(\fcd)$ is a planar set.  

$\partial(\fcd)$ is planar set, and so it is a disjoint collection of circles.  Since $\fcd$ is orientable by the definition of $D^k$-surface, each of its facets are orientable.  By the classification of closed surfaces, it follows that each facet of $\fcd$ is a sphere with punctures or is a $n$-handled torus with punctures, properly embedded in $B^3_+$.  We will now apply ``Corollary 6.2'' from Hempel's text on 3-manifolds, \cite{He}, which says
\begin{center}
\indent \textit{If $F$ is a 2-sided incompressible surface (properly embedded) in a 3-manifold $M$,\\ then ker$\left(\pi_1 (F) \hookrightarrow \pi_1(M)\right)=1$.}\\
\end{center}
Since $\pi_1(B^3_+) \cong 1$, the above result implies that $\pi_1 (\{\textrm{the facets of }\fcd\}) \cong 1$.  There cannot exist facets which are $n$-handled tori nor can there exists facets which are spheres with two or more punctures because this would contradict that $\pi_1 (\{\textrm{the facets of }\fcd\}) \cong 1$.  Therefore, all facets are spheres with a single puncture.  In other words, each facet is a disk. \qed

\begin{lemma}
\label{dk_disk}
Every $D^k$-surface is equal to a linear combination of $D^k$-surfaces whose facets are all disks, and this representation is unique.
\end{lemma}

\emph{Proof.}  If $F$ were a $D^k$-surface for a split link, then it would suffice to show that each connected component of $F$ equaled a linear combination of $D^k$-surfaces whose facets were all disks.  Thus, we may suppose that our link is not a split link.

Let $F$ be a $D^k$-surface for a non-split link $L = \partial(F)$.  If $F$ contains any components without boundary, then the relations (S0), (S1), (D2) and (NC) can be used to remove such components and replace $F$ with a connected $D^k$-surface, $F' = F$ ($F'$ is connected since it contains no components without boundary and $\partial(F') = L$ is connected).  We will proceed by showing that $F'$ can be written as a linear combination of $D^k$-surfaces whose facets are all disks.

By the repeated use of the relation (NC), replace $F'$ with a (finite) linear combination of $D^k$-surfaces which have incompressible facets (this process terminates after a finite number of steps because compressing either increases Euler characteristic or leaves Euler characteristic the same while increasing the number of components).  Use the relations (S0) and (S1) to remove and replace any $D^k$-surfaces in this sum which contain components without boundary.  Hence, we have
$$
F' = \sum_{n=1}^N a_n F_n
$$
where each $a_n \in \Z$ and each $F_n$ is a $D^k$-surface whose facets are not spheres and are incompressible.  By Lemma \ref{facets_are_disks}, each $D^k$-surface $F_n$ in this sum has the property that all of its facets are disks.  Therefore, $F$ can be written as the linear combination of $D^k$-surfaces whose facets are all disks.

It remains to be seen that this representation is unique.  Let $F$ be a $D^k$-surface in $S^3 = B^3_+ \cup B^3_-$ such that $F = \sum_{n=1}^N a_n F_n$ is a decomposition of $F$ into $D^k$-surfaces with disks as facets.  Since the relations (S0), (S1), (D2) and (NC) do not affect a surface near its boundary or its crosscuts, each $D^k$-surface $F_n$ in this linear combination is equal to $F$ inside of a regular neighborhood of the boundary union the crosscuts.  The only thing left to show is that there is only one way to place the facets (disks) of the $F_n$ up to isotopy.

Since a $D^k$-surface and its cross-dual only differ near crosscuts, it suffices to show that there is only one way to attach the facets (disks) of each cross-dual surface to its skeleton.  By the definition of $D^k$-surface, each cross-dual surface can be isotoped so that its skeleton is embedded in $\Sigma = B^3_+ \cap B^3_-$ such that each facet (disk) is in $B^3_+$.  When viewed in this way, we see that attaching the facets of the cross-dual to the skeleton is the same as placing a system of properly embedded disks with predetermined boundary in a closed 3-ball.  Since the process of placing a system of properly embedded disks with predetermined boundary in a closed 3-ball is unique up to isotopy, we are done. \qed
\bigskip

Now we recall the results from Section \ref{6} concerning states surfaces (defined in Definition \ref{statesurface}).  Proposition \ref{state_surf_prop} says that every state surface is a $D^k$-surface.  Given a projection of a link and Khovanov chain complex corresponding to that link projection, Propositions \ref{inj_prop1} and \ref{inj_prop2} prove that there is an injective chain map $\iota$ from the Khovanov chain complex into the diagramless chain complex.  Below we show that every $D^k$-surface can be written as a linear combination of state surfaces, meaning that the diagramless complex contains copies of embedded Khovanov chain complexes and nothing more.  

\begin{prop}
\label{big_prop}
Every $D^k$-surface is equal to a linear combination of state surfaces, and this representation is unique.
\end{prop}

\emph{Proof.}  Let $F$ be a $D^k$-surface.  Using Lemma \ref{dk_disk}, write $F$ as a linear combination of $D^k$-surfaces whose facets are all disks:
$$
F = \sum_{n=1}^N a_n F_n .
$$
Let $F_n$ be an arbitrary $D^k$-surface in this sum.  It suffices to show that $F_n$ is equal to a state surface.  

By the definition of $D^k$-surface there exists an embedded oriented 2-sphere $\Sigma \subseteq S^3$ and embedded 3-balls $B^3_+, B^3_- \subseteq S^3$ such that 
\begin{itemize}
\item [(1)] $\Sigma = B^3_+ \cap B^3_-$, 
\item [(2)] skel$(\fcd_n) \subseteq \Sigma$, 
\item [(3)] all of the locally oriented pieces of surface of $\fcd_n$ agree with the orientation of $\Sigma$, and 
\item [(4)] $\fcd_n - \textrm{skel}(\fcd_n) = \{\textrm{the facets of }\fcd_n\} \subseteq B^3_+$.  
\end{itemize}
We know that the facets of $F_n$ are all disks, so Lemma \ref{facet_lemma} tells us that all of the facets of $\fcd_n$ are disks as well.  Hence, $\fcd_n$ is comprised of properly embedded disks in $B^3_+$ with locally oriented rectangular pieces of surface in $\Sigma$ which are connected to the boundary of two (possibly the same) disks.  Thus $\fcd_n$ can be used as the cross-dual surface in Definition \ref{statesurface} to build a state surface, implying that $F_n$ is a state surface. 

Lemma \ref{dk_disk} tells us that any $D^k$-surface is equal to a unique linear combination of $D^k$-surfaces with facets that are all disks, and we just showed that $D^k$-surfaces with disks as facets are in fact state surfaces.  Therefore any decomposition of a $D^k$-surface into state surfaces must be unique.
\qed
\bigskip

The diagramless chain complex is built from $D^k$-surfaces, so an arbitrary element of the chain complex is a linear combination of $D^k$-surfaces.  By Proposition \ref{big_prop}, we may write a given linear combination of $D^k$-surfaces as a linear combination of state surfaces.  Hence, an arbitrary element of the diagramless chain complex can be written as a linear combination of state surfaces.  

We know that the diagramless complex for a link $L$ consists only of embedded Khovanov complexes for diagrams of the link $L$.  We would like to say that the diagramless complex of a link is equal to the direct sum of some number of embedded Khovanov complexes for that link.  This can be achieved by showing linear independence between certain link diagram equivalence classes.  The appropriate equivalence classes are defined below.

\begin{defn}
\label{e_class}
Let $\Sigma$ be a 2-sphere embedded in $S^3$, and let $D, D'$ be link diagrams in $\Sigma$.  We say that $D$ is equivalent to $D$ and write $D \sim D'$ if there exist enhanced Kauffman states $S_D$ and $S_{D'}$ such that $S_D = S_{D'}$ as $D^k$-surfaces.  Denote the equivalence class of $D$ by $[D]$.  
\end{defn}

\begin{remark}
It should be pointed out that if two link diagrams $D$ and $D'$ differ only by a 2-space isotopy, then $[D] = [D']$.  

What isn't obvious is that the converse does not hold.  For a counterexample, consider a split link diagram $D = D_1 \sqcup D_2$ in the 2-sphere $\Sigma$, where $D_1$ and $D_2$ are sufficiently complicated knot diagrams.  Up to isotopy, there are many different ways to view $D = D_1 \sqcup D_2$ because we could put the knot diagram $D_1$ in any of the different components of $\Sigma - D_2$. 
\end{remark}

\begin{remark}
\label{ec_remark}
By construction, each state surface corresponds to exactly one link diagram equivalence class.  When needed, we denote a such a correspondence by a superscript on the state surface.  For example, a state surface corresponding to the link diagram equivalence class $[D_1]$ could be denoted by $F^{D_1}$ or $F^1$.
\end{remark}

Using these equivalence classes of link diagrams, we are able to show that the diagramless complex breaks as a direct sum of subcomplexes which correspond to distinct link diagram equivalence classes.  This result and the analogous result on the level of homology are stated and proved below.

\begin{prop}
Let $L$ be link and $\mathcal{DC}(L)$ be the diagramless complex of that link, as described in Section \ref{dhsection}.  Denote the distinct link diagram equivalence classes for $L$ by $[D_1],[D_2],...,[D_N]$.  Then we have that 
$$
\mathcal{DC}(L) = \bigoplus_{n=1}^{N} \mathcal{DC}_n(L)
$$
where $\mathcal{DC}_n(L)$ is the subcomplex spanned by the state surfaces which correspond to the link diagram equivalence class $[D_n]$.
\end{prop}

\emph{Proof.}  Consider an arbitrary element $F$ of $\mathcal{DC}(L)$.  Through the use of Proposition \ref{big_prop}, we can write $F$ as the linear combination of state surfaces.  Since each state surface corresponds to exactly one link diagram equivalence class, we may group the linear combination of state surfaces by equivalence class and write
\begin{equation*}
F = \sum_{n=1}^{N} \sum_{m = 1}^{M_n} a_{m}^{n} F^n_m = \left( \sum_{m = 1}^{M_1} a_{m}^{1} F^1_m \right) + \left(  \sum_{m = 1}^{M_2} a_{m}^{2} F^2_m \right) + \cdots + \left(  \sum_{m = 1}^{M_N} a_{m}^{N} F^N_m \right)
\tag{*}
\end{equation*}

Above we have $F$, an arbitrary element of the diagramless chain complex, written as the $N$-term sum of linear combinations of state surfaces corresponding to $N$ distinct link diagram equivalence classes.  To show that $\mathcal{DC}(L)$ breaks as the direct sum of the subcomplexes $\mathcal{DC}_1(L), ..., \mathcal{DC}_N(L)$, it suffices to show that if $F=0$, then each of the $N$ terms in equation (*) equal zero.  That is,we must show that if $F=0$, then $\sum_{m = 1}^{M_n} a_{m}^{N} F^N_m = 0$ for $0 \leq n \leq N$.

Suppose $F=0$.  If each of the $N$ terms in equation (*) equal zero, we are done.  Otherwise, at least two of the $N$ terms in (*) are non-zero, implying that there is cancellation between terms corresponding to different link diagram equivalent classes.  As noted in Remark \ref{ec_remark}, each state surface corresponds to exactly one link diagram equivalence class, so direct cancellation between different link diagram equivalent classes is not possible.  That is, any cancellation that takes place between different link diagram equivalent classes must involve the relations (S1), (S2), (NC) or (D2).  However, this would imply that the use of these relations to represent a $D^k$-surface as a linear combination state surfaces is not unique.  This contradicts Proposition \ref{big_prop}. \qed
\bigskip

As a corollary, we have the following theorem.

\begin{thm}
\label{decomp_thm}
Let $L$ be a link and let $\mathcal{KH}(L)$ denote the Khovanov homology of that link.  The diagramless homology of a link is equal to the direct sum of some number of copies of Khovanov homology for that link.  That is,
$$
\mathcal{DH}_k^i(L) = \bigoplus_{n=1}^N \mathcal{KH}^i(L) = N \cdot \mathcal{KH}^i(L)
$$
where $k$ is the number of crosscuts and $i$ is the homological grading.
\end{thm}

\section{Examples}

In Section \ref{newsection}, the relationship between the diagramless homology and the Khovanov homology of a link was uncovered; the diagramless homology of a link $L$ is equal to the direct sum of some number of copies of the Khovanov homology of $L$.  Since the Khovanov homology of a link is already well-studied and is relatively easy to calculate (see Bar-Natan's calculations in \cite{BN2}), we will only bother with demonstrating a few calculations of the diagramless homology.  In Section ???, the value of $N$, the number of copies of Khovanov homology in the diagramless homology of a given link, is discussed.

\subsection{The Diagramless Homology of the Unknot with \texorpdfstring{$k$}{\textit{k}} = 0 crosscuts}
\label{51}

With no crosscuts, a (connected) $D^k$-surface or a non-split link will only have one facet.  By Lemma \ref{dk_disk}, we know that it suffices to consider $D^k$-surfaces whose facets are disks.  Therefore, the only surface we need to consider is a disk.

Recall that facets of surfaces are allowed to decorated by dots.  Nontrivially, surfaces can only be decorated by zero or one dot per facet due to the (D2) relation.  Therefore, the dotless disk and the once dotted disk are the only two surfaces considered when calculating the $k=0$ diagramless homology of the unknot.  The calculation of the diagramless homology is simple and is show below.

\begin{center}
\begin{tabular}{c}
	\includegraphics[scale=2]{u0.pdf}\\
	\ \\
	\includegraphics[scale=2]{u0.pdf}
	\put(-21,45){\textcolor{blue}{$\bullet$}}
\end{tabular}
\put(-132,15){\vector(2,1){66}}
\put(-132,-5){\vector(2,-1){66}}
\put(3,48){\vector(2,-1){66}}
\put(3,-38){\vector(2,1){66}}
\put(-152,-3){{\Huge 0}}
\put(80,-3){{\Huge 0}}
\put(-179,89){$\overbrace{\ \ \ \ \ \ \ \ \ \ \ \ \ \ \ \ \ \ \ }$}
\put(-64,89){$\overbrace{\ \ \ \ \ \ \ \ \ \ \ \ \ \ \ \ \ \ \ }$}
\put(54,89){$\overbrace{\ \ \ \ \ \ \ \ \ \ \ \ \ \ \ \ \ \ \ }$}
\put(-161,99){$I = -1$}
\put(-45,99){$I = 0$}
\put(74,99){$I = 1$}
\put(-188,-101)
{
	\begin{tabular}{c}
		$\underbrace{\ \ \ \ \ \ \ \ \ \ \ \ \ \ \ \ \ \ \ }$\\
		ker/im\\
		$\shortparallel$\\
		$\langle \emptyset \rangle$
	\end{tabular}
}
\put(-74,-101)
{
	\begin{tabular}{c}
		$\underbrace{\ \ \ \ \ \ \ \ \ \ \ \ \ \ \ \ \ \ \ }$\\
		ker/im\\
		$\shortparallel$\\
		$\langle \1,x \rangle$\\
	\end{tabular}
}
\put(44,-101)
{
	\begin{tabular}{c}
		$\underbrace{\ \ \ \ \ \ \ \ \ \ \ \ \ \ \ \ \ \ \ }$\\
		ker/im\\
		$\shortparallel$\\
		$\langle \emptyset \rangle$\\
	\end{tabular}
}
\end{center}

The notation in the above calculation needs explanation.  Recall that the graded $\Z$-module $M_F := M^{\otimes n}$ is associated to every $D^k$-surface $F$, where $n$ is the number of facets of $F$.  We defined $M$ to be the graded $\Z$-module $\Z[x] / (x^2)$, which is spanned by $\1$ and $x$.  Finally, recall that $\1$ is identified to a facet without a dot and $x$ is identified to a facet with a dot.

Since the dottedness of facets is kept track of by the choice of generator $\1$ or $x$, displaying dotted versions of surfaces can be circumvented.  In light of this fact, we will use the convention that the completely undotted version of a surface will represent all possible dot decorations of that surface.  This will make large diagrams of surfaces much more manageable.  With this convention, the diagram of surfaces depicting the calculation of the unknot for $k=0$ now becomes:

\begin{center}
\includegraphics[scale=2]{u0.pdf}
\put(-125,34){\vector(2,0){66}}
\put(8,34){\vector(2,0){66}}
\put(-147,27){{\Huge 0}}
\put(85,27){{\Huge 0}}
\put(-174,89){$\overbrace{\ \ \ \ \ \ \ \ \ \ \ \ \ \ \ \ \ \ \ }$}
\put(-59,89){$\overbrace{\ \ \ \ \ \ \ \ \ \ \ \ \ \ \ \ \ \ \ }$}
\put(59,89){$\overbrace{\ \ \ \ \ \ \ \ \ \ \ \ \ \ \ \ \ \ \ }$}
\put(-156,99){$I = -1$}
\put(-40,99){$I = 0$}
\put(79,99){$I = 1$}
\put(-183,-34)
{
	\begin{tabular}{c}
		$\underbrace{\ \ \ \ \ \ \ \ \ \ \ \ \ \ \ \ \ \ \ }$\\
		ker/im\\
		$\shortparallel$\\
		$\langle \emptyset \rangle$
	\end{tabular}
}
\put(-69,-34)
{
	\begin{tabular}{c}
		$\underbrace{\ \ \ \ \ \ \ \ \ \ \ \ \ \ \ \ \ \ \ }$\\
		ker/im\\
		$\shortparallel$\\
		$\langle \1,x \rangle$\\
	\end{tabular}
}
\put(49,-34)
{
	\begin{tabular}{c}
		$\underbrace{\ \ \ \ \ \ \ \ \ \ \ \ \ \ \ \ \ \ \ }$\\
		ker/im\\
		$\shortparallel$\\
		$\langle \emptyset \rangle$\\
	\end{tabular}
}
\end{center}

Using the definitions of the four induces (given in Section \ref{3} of this paper), the gradings of the above homology generators $\1$ and $x$ are easily calculated.  The gradings for $\1$ and $x$ are $(I,J,K,B) = (0,-1,0,0)$ and $(I,J,K,B) = (0,1,0,0)$ respectively.  Note that since $M= \Z [x] /(x^2)$ is a $\Z$-module, $\langle \1 \rangle \cong \Z $ and $\langle x \rangle \cong \Z$.

We will use $\D_k^i(L)$ to denote the diagramless homology of the link $L$ present in the gradings $K=k$ and $I=i$.  The notation $\D^{i,j,b}_k(L)$ can be used to distinguish between different $J$ and $B$ gradings when necessary.  With this notation, we can summarize the results for the unknot, $[ \bigcirc ]$, with $k=0$ as follows.

$$
\D_0^i ([ \bigcirc ]) =  \left\{ \begin{array}{cc} \Z \oplus \Z & \textrm{for } \ i=0\\ 0 & \textrm{otherwise} \end{array} \right .
$$

\subsection{The Diagramless Homology of the Unknot with \texorpdfstring{$k$}{\textit{k}} = 1 crosscut}
\label{52}

In this example we will find that the diagramless homology of the unknot with $k=1$ is equal to the direct sum of two copies of Khovanov homology.

By Proposition \ref{big_prop}, it suffices to only determine the state surfaces.  Since state surfaces with one crosscut are built from enhanced Kauffman states that come from diagrams with one crossing, there are a small number of such surfaces to calculate.  

Without the consideration of dots, there are only two types of state surfaces for the unknot with one crosscut: the disk and the Mobius strip.  The disk may contain an \textcolor{darkgreen}{active} or an \textcolor{red}{inactive} crosscut, but the type of Mobius strip (left handed or right handed) determines the type of crosscut allowed.  Also, each facet of a surface is allowed to be decorated with (at most) one dot.  After this consideration, while distinguishing between different facets, we have that there are 12 distinct state surfaces for the unknot (up to isotopy) with one crosscut.  The 12 state surfaces for the unknot with one crosscut are given below with their respective $(I,J,K,B)$ values.

$$
\begin{array}{cccc}
		\ 
	&
		(0,-1,1,1)
	&
		(0,-1,1,0)
	&
		\ 
\\
		\ 
	&
		\ 
		\includegraphics[scale=1.5]{u0.pdf}
		\thicklines
		\color{darkgreen}
		\put(-37.5,22){\vector(3,1){35.3}}
		\color{black}
		\ 
	&
		\ 
		\includegraphics[scale=1.5]{u0.pdf}
		\thicklines
		\color{red}
		\put(-37.5,22){\vector(3,1){35.3}}
		\color{black}
		\ 
	&
		\ 
\\
		\ 
	&
	 	\ 
	&
		\ 
	&
		\ 
\\
		(-1,1,1,0)
	&
		(0,1,1,1)
	&
		(0,1,1,0)
	&
		(1,-1,1,1)
\\
		\reflectbox{ \includegraphics[scale=1.5]{u1.pdf} }
		\thicklines
		\color{darkgreen}
		\put(-26.5,39){\vector(1,1){12.3}}
		\color{black}
	&
		\ 
		\includegraphics[scale=1.5]{u0.pdf}
		\thicklines
		\color{darkgreen}
		\put(-37.5,22){\vector(3,1){35.3}}
		\color{blue}
		\put(-20,10){$\bullet$}
		\color{black}
		\ 
	&
		\ 
		\includegraphics[scale=1.5]{u0.pdf}
		\thicklines
		\color{red}
		\put(-37.5,22){\vector(3,1){35.3}}
		\color{blue}
		\put(-20,10){$\bullet$}
		\color{black}
		\ 
	&
		\ 
		\includegraphics[scale=1.5]{u1.pdf} 
		\thicklines
		\color{red}
		\put(-23.2,39){\vector(1,1){12.3}}
		\color{black}
		\ 
\\
		\ 
	&
	 	\ 
	&
		\ 
	&
		\ 
\\
		(-1,3,1,0)
	&
		(0,1,1,1)
	&
		(0,1,1,0)
	&
		(1,1,1,1)
\\
		\reflectbox{ \includegraphics[scale=1.5]{u1.pdf} }
		\thicklines
		\color{darkgreen}
		\put(-26.5,39){\vector(1,1){12.3}}
		\color{blue}
		\put(-20,20){$\bullet$}
		\color{black}
	&
		\ 
		\includegraphics[scale=1.5]{u0.pdf}
		\thicklines
		\color{darkgreen}
		\put(-37.5,22){\vector(3,1){35.3}}
		\color{blue}
		\put(-25,37){$\bullet$}
		\color{black}
		\ 
	&
		\ 
		\includegraphics[scale=1.5]{u0.pdf}
		\thicklines
		\color{red}
		\put(-37.5,22){\vector(3,1){35.3}}
		\color{blue}
		\put(-25,37){$\bullet$}
		\color{black}
		\ 
	&
		\ 
		\includegraphics[scale=1.5]{u1.pdf} 
		\thicklines
		\color{red}
		\put(-23.2,39){\vector(1,1){12.3}}
		\color{blue}
		\put(-16,20){$\bullet$}
		\color{black}
		\ 
\\
		\ 
	&
	 	\ 
	&
		\ 
	&
		\ 
\\
		\ 
	&
		(0,3,1,1)
	&
		(0,3,1,0)
	&
		\ 
\\
		\ 
	&
		\ 
		\includegraphics[scale=1.5]{u0.pdf}
		\thicklines
		\color{darkgreen}
		\put(-37.5,22){\vector(3,1){35.3}}
		\color{blue}
		\put(-25,37){$\bullet$}
		\put(-20,10){$\bullet$}
		\color{black}
		\ 
	&
		\ 
		\includegraphics[scale=1.5]{u0.pdf}
		\thicklines
		\color{red}
		\put(-37.5,22){\vector(3,1){35.3}}
		\color{blue}
		\put(-25,37){$\bullet$}
		\put(-20,10){$\bullet$}
		\color{black}
		\ 
	&
		\ 
\end{array}
$$

Typically, the visual presence of dots will be suppressed and instead we denote a dotted facet by an $x$ and an undotted facet by a $\1$ (this convention was explained at the end of Section \ref{51}).  The $n$ facets of a state surface are labeled by $\underline{1}, ..., \underline{n}$.  An $x$ (resp.\thinspace $\1$) in the $\ell^{\textrm{th}}$ coordinate of a tensor product denotes the presence (resp.\thinspace absence) of a dot on facet $\underline{\ell}$ of the surface.  

With this notation, the 12 different state surfaces with one crosscut can be represented as follows:
$$
\begin{array}{cccc}
		\reflectbox{ \includegraphics[scale=1.5]{u1.pdf} }
		\thicklines
		\color{darkgreen}
		\put(-26.5,39){\vector(1,1){12.3}}
		\color{black}
		\put(-39,10){{\large $\underline{1}$}}
	&
		\ 
		\includegraphics[scale=1.5]{u0.pdf}
		\thicklines
		\color{darkgreen}
		\put(-37.5,22){\vector(3,1){35.3}}
		\color{black}
		\put(-24,38){{\large $\underline{1}$}}
		\put(-22,11){{\large $\underline{2}$}}
		\ 
	&
		\ 
		\includegraphics[scale=1.5]{u0.pdf}
		\thicklines
		\color{red}
		\put(-37.5,22){\vector(3,1){35.3}}
		\color{black}
		\put(-24,38){{\large $\underline{1}$}}
		\put(-22,11){{\large $\underline{2}$}}
		\ 
	&
		\ 
		\includegraphics[scale=1.5]{u1.pdf} 
		\thicklines
		\color{red}
		\put(-23.2,39){\vector(1,1){12.3}}
		\color{black}
		\put(-36,10){{\large $\underline{1}$}}
		\ 
\\
			\ 
	&
		\ \ \1 \otimes \1 \sim  (0,-1,1,1)
	&
		\ \ \1 \otimes \1 \sim  (0,-1,1,0)
	&
		\ 
\\
		\1 \sim  (-1,1,1,0)
	&
		\1 \otimes x \sim  (0,1,1,1)
	&
		\1 \otimes x \sim  (0,1,1,0)
	&
		\ \ \1 \sim  (1,-1,1,1)
\\
		x \sim  (-1,3,1,0)
	&
		x \otimes \1 \sim  (0,1,1,1)
	&
		x \otimes \1 \sim  (0,1,1,0)
	&
		x \sim  (1,1,1,1)
\\
		\ 
	&
		x \otimes x \sim  (0,3,1,1)
	&
		x \otimes x \sim  (0,3,1,0)
	&
		\ 
\end{array}
$$

Now that the state surfaces are known, the diagramless homology for the unknot with $k=1$ can be calculated.  A diagram corresponding to this calculation is given below.

\begin{center}
\begin{tabular}{c}
\ \\
		\reflectbox{ \includegraphics[scale=1.5]{u1.pdf} }
		\thicklines
		\color{darkgreen}
		\put(-29.7,39){\vector(1,1){12.3}}
		\color{black}
		\put(-39,10){{\large $\underline{1}$}}
\end{tabular}
\ \ \ \ \ \ \ \ \ \ \ \ 
\begin{tabular}{c}
\ \\
		\includegraphics[scale=1.5]{u0.pdf}
		\thicklines
		\color{red}
		\put(-40.7,22){\vector(3,1){35.3}}
		\color{black}
		\put(-24,38){{\large $\underline{1}$}}
		\put(-22,11){{\large $\underline{2}$}}
\ \\
\ \\
\ \\
\ \\
		\includegraphics[scale=1.5]{u0.pdf}
		\thicklines
		\color{darkgreen}
		\put(-40.7,22){\vector(3,1){35.3}}
		\color{black}
		\put(-24,38){{\large $\underline{1}$}}
		\put(-22,11){{\large $\underline{2}$}}
\end{tabular}
\put(-100,15){\vector(2,1){38}}
\put(4,-38){\vector(2,1){41}}
\put(-177,82){$\overbrace{\ \ \ \ \ \ \ \ \ \ \ \ \ \ \ \ \ \ \ }$}
\put(-61,82){$\overbrace{\ \ \ \ \ \ \ \ \ \ \ \ \ \ \ \ \ \ \ }$}
\put(59,82){$\overbrace{\ \ \ \ \ \ \ \ \ \ \ \ \ \ \ \ \ \ \ }$}
\put(-159,92){$I = -1$}
\put(-42,92){$I = 0$}
\put(79,92){$I = 1$}
\put(-187,-116)
{
	\begin{tabular}{c}
		$\underbrace{\ \ \ \ \ \ \ \ \ \ \ \ \ \ \ \ \ \ \ }$\\
		ker/im\\
		$\shortparallel$\\
		$\langle \emptyset \rangle / \langle \emptyset \rangle$\\
		$\shortparallel$\\
		$\langle \emptyset \rangle$\\
	\end{tabular}
}
\put(-126,-118)
{
	\begin{tabular}{c}
		$\ \ \underbrace{\ \ \ \ \ \ \ \ \ \ \ \ \ \ \ \ \ \ \ }$\\
		\ \ ker/im\\
		$\ \ \shortparallel$\\
		$\frac{\langle \1 \otimes x, \1 \otimes \1, x \otimes x, x \otimes \1 \rangle}{\langle \1 \otimes x + x \otimes \1, x \otimes x \rangle} \bigoplus \frac{\langle \1 \otimes x -  x \otimes \1, x \otimes x \rangle}{\langle \emptyset \rangle}$\\
		$\ \ \shortparallel$\\
		$\langle \1 \otimes x,\1 \otimes \1 \rangle \oplus \langle \1 \otimes x -  x \otimes \1, x \otimes x \rangle$\\
	\end{tabular}
}
\put(48,-116)
{
	\begin{tabular}{c}
		$\underbrace{\ \ \ \ \ \ \ \ \ \ \ \ \ \ \ \ \ \ \ }$\\
		ker/im\\
		$\shortparallel$\\
		$\langle \1, x \rangle / \langle \1, x \rangle$\\
		$\shortparallel$\\
		$\langle \emptyset \rangle$\\
	\end{tabular}
}
\ \ \ \ \ \ \ \ \ \ \ \ \ \ \ 
\begin{tabular}{c}
\ \\
		\includegraphics[scale=1.5]{u1.pdf}
		\thicklines
		\color{red}
		\put(-26.4,39){\vector(1,1){12.3}}
		\color{black}
		\put(-36,10){{\large $\underline{1}$}}
\end{tabular}
\end{center}

In the above diagram we see how the diagramless chain complex splits into a direct sum of Khovanov chain complexes, each with a different value of the index $B$.  For $B=0$ (the top complex), the only nontrivial homology is in the $I=0$ grading.  Here $\langle \1 \otimes x,\1 \otimes \1 \rangle \cong \Z \oplus \Z$.  For $B=1$ (the bottom complex), the only nontrivial homology is again in the $I=0$ grading.  Here $\langle \1 \otimes x -  x \otimes \1, x \otimes x \rangle \cong \Z \oplus \Z$ as well.

Using the notation introduced at the end of Section \ref{51}, the diagramless homology for the unknot with $k=1$ is summarized below.

$$
\D^{i,b}_1 ([ \bigcirc ]) \cong  \left\{ \begin{array}{cc} \Z \oplus \Z & \textrm{for } \ i=0, \ b=0\\ \Z \oplus \Z & \textrm{for } \ i=0, \ b=1\\ 0 & \textrm{otherwise} \end{array} \right .
$$

\subsection{The Diagramless Homology of the Unknot with \texorpdfstring{$k$}{\textit{k}} = 2 crosscuts}
\label{53}

Just as in the previous section, realize that it suffices to only consider states surface by Proposition \ref{big_prop}.  Hence, we only need to find the enhanced Kauffman states that come from diagrams for the unknot with 2 crossings.  In this example, there are many state surfaces.  Below we list the undotted version of each surface.

$$
\begin{array}{|cc|c||c|c|} \hline
		\ \ \ \ \ \
		\put(-13,23){ \begin{tabular}{l}
			$I=$ \textbf{-}2 \\
			$J=$ 3 \\
			$K=$ 2 \\
			$B=$ 0
			\end{tabular} }
	& \ \ \ \ 
	&
		\reflectbox{ \includegraphics[scale=1.5]{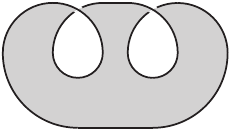} }
		\thicklines
		\color{darkgreen}
		\put(-26,39){\vector(1,1){12.3}}
		\put(-80,39){\vector(-1,1){12.3}}
		\color{black}
		\put(-2,25){$\cong$}
		\ 
		\reflectbox{ \includegraphics[scale=1.5]{u2.pdf} }
		\thicklines
		\color{darkgreen}
		\put(-14,51){\vector(-1,-1){12.3}}
		\put(-92,51){\vector(1,-1){12.3}}
		\color{black}
	&
		\ 
			\put(-13,23){ \begin{tabular}{l}
			$I=$ \textbf{-}2 \\
			$J=$ 1 \\
			$K=$ 2 \\
			$B=$ 2
			\end{tabular} }
		\ \ \ \ \ \ \ 
	&
		\includegraphics[scale=1.5]{u0.pdf}
		\thicklines
		\color{darkgreen}
		\put(-37.2,17){\vector(1,0){34.5}}
		\put(-37.2,37){\vector(1,0){34.5}}
		\color{black}
		\put(1,25){$\cong$}
		\ \ \ 
		\includegraphics[scale=1.5]{u0.pdf}
		\thicklines
		\color{darkgreen}
		\put(-2.7,17){\vector(-1,0){34.5}}
		\put(-2.7,37){\vector(-1,0){34.5}}
		\color{black}
\\ \hline
		\ \ \ \ \ \
		\put(-13,23){ \begin{tabular}{l}
			$I=$ \textbf{-}2 \\
			$J=$ 3 \\
			$K=$ 2 \\
			$B=$ 0
			\end{tabular} }
	& \ \ \ \ 
	&
		\reflectbox{ \includegraphics[scale=1.5]{u2.pdf} }
		\thicklines
		\color{darkgreen}
		\put(-14,51){\vector(-1,-1){12.3}}
		\put(-80,39){\vector(-1,1){12.3}}
		\color{black}
		\put(-2,25){$\cong$}
		\ 
		\reflectbox{ \includegraphics[scale=1.5]{u2.pdf} }
		\thicklines
		\color{darkgreen}
		\put(-26,39){\vector(1,1){12.3}}
		\put(-92,51){\vector(1,-1){12.3}}
		\color{black}
	&
		\ 
			\put(-13,23){ \begin{tabular}{l}
			$I=$ \textbf{-}2 \\
			$J=$ 1 \\
			$K=$ 2 \\
			$B=$ 2
			\end{tabular} }
		\ \ \ \ \ \ \ 
	&
		\includegraphics[scale=1.5]{u0.pdf}
		\thicklines
		\color{darkgreen}
		\put(-37.2,17){\vector(1,0){34.5}}
		\put(-2.7,37){\vector(-1,0){34.5}}
		\color{black}
		\put(1,25){$\cong$}
		\ \ \ 
		\includegraphics[scale=1.5]{u0.pdf}
		\thicklines
		\color{darkgreen}
		\put(-2.7,17){\vector(-1,0){34.5}}
		\put(-37.2,37){\vector(1,0){34.5}}
		\color{black}
\\ \hline
		\ \ \ \ \ \
		\put(-13,23){ \begin{tabular}{l}
			$I=$ 0 \\
			$J=$ 1 \\
			$K=$ 2 \\
			$B=$ 1
			\end{tabular} }
	& \ \ \ \ 
	&
		\includegraphics[scale=1.5]{u00.pdf}
		\thicklines
		\color{darkgreen}
		\put(-11,51){\vector(-1,-1){12.3}}
		\color{red}
		\put(-89,51){\vector(1,-1){12.3}}
		\color{black}
		\put(1,25){$\cong$}
		\ \ \ 
		\includegraphics[scale=1.5]{u00.pdf}
		\thicklines
		\color{darkgreen}
		\put(-23,39){\vector(1,1){12.3}}
		\color{red}
		\put(-77,39){\vector(-1,1){12.3}}
		\color{black}
	&
		\ 
			\put(-13,23){ \begin{tabular}{l}
			$I=$ 0 \\
			$J=$ \textbf{-}1 \\
			$K=$ 2 \\
			$B=$ 1
			\end{tabular} }
		\ \ \ \ \ \ \ 
	&
		\includegraphics[scale=1.5]{u0.pdf}
		\thicklines
		\color{darkgreen}
		\put(-37.2,17){\vector(1,0){34.5}}
		\color{red}
		\put(-37.2,37){\vector(1,0){34.5}}
		\color{black}
		\put(1,25){$\cong$}
		\ \ \ 
		\includegraphics[scale=1.5]{u0.pdf}
		\thicklines
		\color{darkgreen}
		\put(-2.7,17){\vector(-1,0){34.5}}
		\color{red}
		\put(-2.7,37){\vector(-1,0){34.5}}
		\color{black}
\\ \hline
		\ \ \ \ \ \
		\put(-13,23){ \begin{tabular}{l}
			$I=$ 0 \\
			$J=$ 1 \\
			$K=$ 2 \\
			$B=$ 1
			\end{tabular} }
	& \ \ \ \  &
		\includegraphics[scale=1.5]{u00.pdf}
		\thicklines
		\color{darkgreen}
		\put(-23,39){\vector(1,1){12.3}}
		\color{red}
		\put(-89,51){\vector(1,-1){12.3}}
		\color{black}
		\put(1,25){$\cong$}
		\ \ \ 
		\includegraphics[scale=1.5]{u00.pdf}
		\thicklines
		\color{darkgreen}
		\put(-11,51){\vector(-1,-1){12.3}}
		\color{red}
		\put(-77,39){\vector(-1,1){12.3}}
		\color{black}
	&
		\ 
			\put(-13,23){ \begin{tabular}{l}
			$I=$ 0 \\
			$J=$ \textbf{-}1 \\
			$K=$ 2 \\
			$B=$ 1
			\end{tabular} }
		\ \ \ \ \ \ \ 
	&
		\includegraphics[scale=1.5]{u0.pdf}
		\thicklines
		\color{darkgreen}
		\put(-37.2,17){\vector(1,0){34.5}}
		\color{red}
		\put(-2.7,37){\vector(-1,0){34.5}}
		\color{black}
		\put(1,25){$\cong$}
		\ \ \ 
		\includegraphics[scale=1.5]{u0.pdf}
		\thicklines
		\color{darkgreen}
		\put(-2.7,17){\vector(-1,0){34.5}}
		\color{red}
		\put(-37.2,37){\vector(1,0){34.5}}
		\color{black}
\\ \hline
		\ \ \ \ \ \
		\put(-13,23){ \begin{tabular}{l}
			$I=$ 2 \\
			$J=$ \textbf{-}1 \\
			$K=$ 2 \\
			$B=$ 2
			\end{tabular} }
	& \ \ \ \  &
		\includegraphics[scale=1.5]{u2.pdf}
		\thicklines
		\color{red}
		\put(-23,39){\vector(1,1){12.3}}
		\put(-77,39){\vector(-1,1){12.3}}
		\color{black}
		\put(1,25){$\cong$}
		\ \ \ 
		\includegraphics[scale=1.5]{u2.pdf}
		\thicklines
		\color{red}
		\put(-11,51){\vector(-1,-1){12.3}}
		\put(-89,51){\vector(1,-1){12.3}}
		\color{black}
	&
		\ 
			\put(-13,23){ \begin{tabular}{l}
			$I=$ 0 \\
			$J=$ \textbf{-}1 \\
			$K=$ 2 \\
			$B=$ 0
			\end{tabular} }
		\ \ \ \ \ \ \ 
	&
		\includegraphics[scale=1.5]{u0.pdf}
		\thicklines
		\color{red}
		\put(-37.2,17){\vector(1,0){34.5}}
		\put(-37.2,37){\vector(1,0){34.5}}
		\color{black}
		\put(1,25){$\cong$}
		\ \ \ 
		\includegraphics[scale=1.5]{u0.pdf}
		\thicklines
		\color{red}
		\put(-2.7,17){\vector(-1,0){34.5}}
		\put(-2.7,37){\vector(-1,0){34.5}}
		\color{black}
\\ \hline
		\ \ \ \ \ \
		\put(-13,23){ \begin{tabular}{l}
			$I=$ 2 \\
			$J=$ \textbf{-}1 \\
			$K=$ 2 \\
			$B=$ 2
			\end{tabular} }
	& \ \ \ \  &
		\includegraphics[scale=1.5]{u2.pdf}
		\thicklines
		\color{red}
		\put(-11,51){\vector(-1,-1){12.3}}
		\put(-77,39){\vector(-1,1){12.3}}
		\color{black}
		\put(1,25){$\cong$}
		\ \ \ 
		\includegraphics[scale=1.5]{u2.pdf}
		\thicklines
		\color{red}
		\put(-23,39){\vector(1,1){12.3}}
		\put(-89,51){\vector(1,-1){12.3}}
		\color{black}
	&
		\ 
			\put(-13,23){ \begin{tabular}{l}
			$I=$ 0 \\
			$J=$ \textbf{-}1 \\
			$K=$ 2 \\
			$B=$ 0
			\end{tabular} }
		\ \ \ \ \ \ \ 
	&
		\includegraphics[scale=1.5]{u0.pdf}
		\thicklines
		\color{red}
		\put(-37.2,17){\vector(1,0){34.5}}
		\put(-2.7,37){\vector(-1,0){34.5}}
		\color{black}
		\put(1,25){$\cong$}
		\ \ \ 
		\includegraphics[scale=1.5]{u0.pdf}
		\thicklines
		\color{red}
		\put(-2.7,17){\vector(-1,0){34.5}}
		\put(-37.2,37){\vector(1,0){34.5}}
		\color{black}
\\ \hline
		\ \ \ \ \ \
		\put(-13,23){ \begin{tabular}{l}
			$I=$ 1 \\
			$J=$ \textbf{-}1 \\
			$K=$ 2 \\
			$B=$ 1
			\end{tabular} }
	& \ \ \ \  &
		\includegraphics[scale=1.5]{u1.pdf} 
		\thicklines
		\color{red}
		\put(-11,51){\vector(-1,-1){12.3}}		
		\put(-59.5,14){\vector(1,0){52}}
		\color{black}
		\put(1,25){$\cong$}
		\ \ \ 
		\includegraphics[scale=1.5]{u1.pdf} 
		\thicklines
		\color{red}
		\put(-23.2,39){\vector(1,1){12.3}}
		\put(-59.5,14){\vector(1,0){52}}
		\color{black}
	&
		\ 
			\put(-13,23){ \begin{tabular}{l}
			$I=$ 1 \\
			$J=$ \textbf{-}1 \\
			$K=$ 2 \\
			$B=$ 1
			\end{tabular} }
		\ \ \ \ \ \ \ 
	&
		\reflectbox{ \includegraphics[scale=1.5]{u1.pdf} }
		\thicklines
		\color{darkgreen}
		\put(-14,51){\vector(-1,-1){12.3}}		
		\put(-62.5,14){\vector(1,0){52}}
		\color{black}
		\put(-2,25){$\cong$}
		\ 
		\reflectbox{ \includegraphics[scale=1.5]{u1.pdf} }
		\thicklines
		\color{darkgreen}
		\put(-26.2,39){\vector(1,1){12.3}}
		\put(-62.5,14){\vector(1,0){52}}
		\color{black}
\\ \hline
		\ \ \ \ \ \
		\put(-13,23){ \begin{tabular}{l}
			$I=$ 1 \\
			$J=$ \textbf{-}1 \\
			$K=$ 2 \\
			$B=$ 1
			\end{tabular} }
	& \ \ \ \  &
		\includegraphics[scale=1.5]{u1.pdf} 
		\thicklines
		\color{red}
		\put(-11,51){\vector(-1,-1){12.3}}		
		\put(-7.5,14){\vector(-1,0){52}}
		\color{black}
		\put(1,25){$\cong$}
		\ \ \ 
		\includegraphics[scale=1.5]{u1.pdf} 
		\thicklines
		\color{red}
		\put(-23.2,39){\vector(1,1){12.3}}
		\put(-7.5,14){\vector(-1,0){52}}
		\color{black}
	&
		\ 
			\put(-13,23){ \begin{tabular}{l}
			$I=$ 1 \\
			$J=$ \textbf{-}1 \\
			$K=$ 2 \\
			$B=$ 1
			\end{tabular} }
		\ \ \ \ \ \ \ 
	&
		\reflectbox{ \includegraphics[scale=1.5]{u1.pdf} }
		\thicklines
		\color{darkgreen}
		\put(-14,51){\vector(-1,-1){12.3}}		
		\put(-10.5,14){\vector(-1,0){52}}
		\color{black}
		\put(-2,25){$\cong$}
		\ 
		\reflectbox{ \includegraphics[scale=1.5]{u1.pdf} }
		\thicklines
		\color{darkgreen}
		\put(-26.2,39){\vector(1,1){12.3}}
		\put(-10.5,14){\vector(-1,0){52}}
		\color{black}
\\ \hline
		\ \ \ \ \ \
		\put(-13,23){ \begin{tabular}{l}
			$I=$ 1 \\
			$J=$ \textbf{-}1 \\
			$K=$ 2 \\
			$B=$ 2
			\end{tabular} }
	& \ \ \ \  &
		\includegraphics[scale=1.5]{u1.pdf} 
		\thicklines
		\color{red}
		\put(-23.2,39){\vector(1,1){12.3}}
		\color{darkgreen}
		\put(-7.5,14){\vector(-1,0){52}}
		\color{black}
		\put(1,25){$\cong$}
		\ \ \ 
		\includegraphics[scale=1.5]{u1.pdf} 
		\thicklines
		\color{red}
		\put(-11,51){\vector(-1,-1){12.3}}
		\color{darkgreen}
		\put(-7.5,14){\vector(-1,0){52}}
		\color{black}
	&
		\ 
			\put(-13,23){ \begin{tabular}{l}
			$I=$ \textbf{-}1 \\
			$J=$ 1 \\
			$K=$ 2 \\
			$B=$ 0
			\end{tabular} }
		\ \ \ \ \ \ \ 
	&
		\reflectbox{ \includegraphics[scale=1.5]{u1.pdf} }
		\thicklines
		\color{darkgreen}
		\put(-26.2,39){\vector(1,1){12.3}}
		\color{red}
		\put(-10.5,14){\vector(-1,0){52}}
		\color{black}
		\put(-2,25){$\cong$}
		\ 
		\reflectbox{ \includegraphics[scale=1.5]{u1.pdf} }
		\thicklines
		\color{darkgreen}
		\put(-14,51){\vector(-1,-1){12.3}}
		\color{red}
		\put(-10.5,14){\vector(-1,0){52}}
		\color{black}
\\ \hline
		\ \ \ \ \ \
		\put(-13,23){ \begin{tabular}{l}
			$I=$ 1 \\
			$J=$ \textbf{-}1 \\
			$K=$ 2 \\
			$B=$ 2
			\end{tabular} }
	& \ \ \ \  &
		\includegraphics[scale=1.5]{u1.pdf} 
		\thicklines
		\color{red}
		\put(-23.2,39){\vector(1,1){12.3}}
		\color{darkgreen}
		\put(-59.5,14){\vector(1,0){52}}
		\color{black}
		\put(1,25){$\cong$}
		\ \ \ 
		\includegraphics[scale=1.5]{u1.pdf} 
		\thicklines
		\color{red}
		\put(-11,51){\vector(-1,-1){12.3}}
		\color{darkgreen}
		\put(-59.5,14){\vector(1,0){52}}
		\color{black}
	&
		\ 
			\put(-13,23){ \begin{tabular}{l}
			$I=$ \textbf{-}1 \\
			$J=$ 1 \\
			$K=$ 2 \\
			$B=$ 0
			\end{tabular} }
		\ \ \ \ \ \ \ 
	&
		\reflectbox{ \includegraphics[scale=1.5]{u1.pdf} }
		\thicklines
		\color{darkgreen}
		\put(-26.2,39){\vector(1,1){12.3}}
		\color{red}
		\put(-62.5,14){\vector(1,0){52}}
		\color{black}
		\put(-2,25){$\cong$}
		\ 
		\reflectbox{ \includegraphics[scale=1.5]{u1.pdf} }
		\thicklines
		\color{darkgreen}
		\put(-14,51){\vector(-1,-1){12.3}}
		\color{red}
		\put(-62.5,14){\vector(1,0){52}}
		\color{black}
\\ \hline
\end{array}
$$

The diagramless homology of the unknot for $k=2$ will be calculated using the same conventions as the previous example -- the $n$ facets of a $D^k$-surface are labeled by $\underline{1}, ..., \underline{n}$, and an $x$ (resp.\thinspace $\1$) in the $\ell^{\textrm{th}}$ coordinate of a tensor product denotes the presence (resp.\thinspace absence) of a dot on facet $\underline{\ell}$ of the surface.  For the $k=1$ diagramless homology of the unknot, we saw that the chain complex could be viewed as two subcomplexes separated by different $B$-grading values.  For $k=2$, there are three non-trivial subcomplexes separated by different $B$-grading values, each of which can be further decomposed into two copies of embedded Khovanov complexes.

\bigskip

The $B = 0$ subcomplex:\footnote{Recall that the differential $d$ is defined a weighted sum of the $d_c$ maps, where the weight is either +1 or -1.  If a map $d_c$ has a negative weight, we will note this by putting a small circle at the initial point of the arrow representing that map in our homology diagram.  This is the same convention that Bar-Natan uses in his paper on Khovanov homology, \cite{BN1}. }

\begin{center}
\begin{tabular}{c}
\ \\
		\reflectbox{ \includegraphics[scale=1.5]{u2.pdf} }
		\thicklines
		\color{darkgreen}
		\put(-39.6,34){$c_2$}
		\put(-81.5,34){$c_1$}
		\put(-29,39){\vector(1,1){12.3}}
		\put(-83,39){\vector(-1,1){12.3}}
		\color{black}
		\put(-59,10){{$\underline{1}$}}
\end{tabular}
\ \ \ \ \ \ \ \ 
\begin{tabular}{c}
\ \\
		\reflectbox{ \includegraphics[scale=1.5]{u1.pdf} }
		\thicklines
		\color{darkgreen}
		\put(-29.2,39){\vector(1,1){12.3}}
		\put(-40,34){$c_2$}
		\color{red}
		\put(-74.5,9){$c_1$}
		\put(-14,14){\vector(-1,0){52}}
		\color{black}
		\put(-63,32){{$\underline{1}$}}
		\put(-42,5){{$\underline{2}$}}
\ \\
\ \\
\ \\
\ \\
		\reflectbox{ \includegraphics[scale=1.5]{u1.pdf} }
		\thicklines
		\color{darkgreen}
		\put(-29.2,39){\vector(1,1){12.3}}
		\put(-40.5,34.5){$c_1$}
		\color{red}
		\put(-75,9){$c_2$}
		\put(-14,14){\vector(-1,0){52}}
		\color{black}
		\put(-63,32){{$\underline{1}$}}
		\put(-42,5){{$\underline{2}$}}
\end{tabular}
\put(-125,16){\vector(2,1){38}}
\put(-125,-16){\vector(2,-1){38}}
\put(-3,-37){\vector(2,1){50}}
\put(-3,37){\vector(2,-1){50}}
\put(-6.5,35.5){$\circ$}
\put(-213,82){$\overbrace{\ \ \ \ \ \ \ \ \ \ \ \ \ \ \ \ \ \ \ }$}
\put(-78,82){$\overbrace{\ \ \ \ \ \ \ \ \ \ \ \ \ \ \ \ \ \ \ }$}
\put(44,82){$\overbrace{\ \ \ \ \ \ \ \ \ \ \ \ \ \ \ \ \ \ \ }$}
\put(-195,92){$I = -2$}
\put(-60,92){$I = -1$}
\put(64,92){$I = 0$}
\ \ \ \ \ \ \ \ \ \ \ \ \ \ \ 
\begin{tabular}{c}
\ \\
		\includegraphics[scale=1.5]{u0.pdf}
		\thicklines
		\color{red}
		\put(-40.2,17){\vector(1,0){34.5}}
		\put(-40.2,37){\vector(1,0){34.5}}
		\put(-49.7,13.5){$c_1$}
		\put(-51.2,35){$c_2$}
		\color{black}
		\put(-26,43){{$\underline{1}$}}
		\put(-26,24){{$\underline{2}$}}
		\put(-26,6){{$\underline{3}$}}
\end{tabular} \ \ \ \ \ \ \ \ 
\end{center}

\begin{center}
\begin{tabular}{c}
\ \\
		\reflectbox{ \includegraphics[scale=1.5]{u2.pdf} }
		\thicklines
		\color{darkgreen}
		\put(-17,51){\vector(-1,-1){12.3}}
		\put(-83,39){\vector(-1,1){12.3}}
		\put(-40.4,34){$c_2$}
		\put(-81,34){$c_1$}
		\color{black}
		\put(-59,10){{$\underline{1}$}}
\end{tabular}
\ \ \ \ \ \ \ \ 
\begin{tabular}{c}
\ \\
		\reflectbox{ \includegraphics[scale=1.5]{u1.pdf} }
		\thicklines
		\color{darkgreen}
		\put(-29.2,39){\vector(1,1){12.3}}
		\put(-40.5,34){$c_2$}
		\color{red}
		\put(-66,14){\vector(1,0){52}}
		\put(-75,9.5){$c_1$}
		\color{black}
		\put(-63,32){{$\underline{1}$}}
		\put(-42,5){{$\underline{2}$}}
\ \\
\ \\
\ \\
\ \\
		\reflectbox{ \includegraphics[scale=1.5]{u1.pdf} }
		\thicklines
		\color{darkgreen}
		\put(-29.2,39){\vector(1,1){12.3}}
		\put(-40,34){$c_1$}
		\color{red}
		\put(-75.5,9.5){$c_2$}
		\put(-66,14){\vector(1,0){52}}
		\color{black}		
		\put(-63,32){{$\underline{1}$}}
		\put(-42,5){{$\underline{2}$}}
\end{tabular}
\put(-125,16){\vector(2,1){38}}
\put(-125,-16){\vector(2,-1){38}}
\put(-3,-37){\vector(2,1){50}}
\put(-3,37){\vector(2,-1){50}}
\put(-6.5,35.5){$\circ$}
\put(-223,-132)
{
	\begin{tabular}{c}
		$\underbrace{\ \ \ \ \ \ \ \ \ \ \ \ \ \ \ \ \ \ \ }$\\
		ker/im\\
		$\shortparallel$\\
		$\vdots$\\
		$\shortparallel$\\
		$\langle \emptyset \rangle$\\
		$\oplus$\\
	 	$\langle \emptyset \rangle$
	\end{tabular}
}
\put(-94,-132)
{
	\begin{tabular}{c}
		$\ \ \underbrace{\ \ \ \ \ \ \ \ \ \ \ \ \ \ \ \ \ \ \ }$\\
		\ \ ker/im\\
		$\ \ \shortparallel$\\		
		\ \ $\vdots$\\
		$\ \ \shortparallel$\\
		\ \ $\langle \emptyset \rangle$\\
		\ \ $\oplus$\\
	 	\ \ $\langle \emptyset \rangle$
	\end{tabular}
}
\put(18,-132)
{
	\begin{tabular}{c}
		$\underbrace{\ \ \ \ \ \ \ \ \ \ \ \ \ \ \ \ \ \ \ }$\\
		ker/im\\
		$\shortparallel$\\
		$\vdots$\\
		$\shortparallel$\\
		$\langle \1 \otimes \1 \otimes \1, \ \1 \otimes \1 \otimes x \rangle$\\
		$\oplus$\\
		$\langle \1 \otimes \1 \otimes \1, \ \1 \otimes \1 \otimes x \rangle$
	\end{tabular}
}
\ \ \ \ \ \ \ \ \ \ \ \ \ \ \ 
\begin{tabular}{c}
\ \\
		\includegraphics[scale=1.5]{u0.pdf}
		\thicklines
		\color{red}
		\put(-40.2,17){\vector(1,0){34.5}}
		\put(-5.7,37){\vector(-1,0){34.5}}
		\put(-49.5,13.5){$c_1$}
		\put(-51.2,35.3){$c_2$}
		\color{black}
		\put(-26,43){{$\underline{1}$}}
		\put(-26,24){{$\underline{2}$}}
		\put(-26,6){{$\underline{3}$}}
\end{tabular} \ \ \ \ \ \ \ \ 
\end{center}

\newpage
The $B = 1$ subcomplex:

\begin{center}
\ \ \ \ \ \ \begin{tabular}{c}
\ \\
		\reflectbox{ \includegraphics[scale=1.5]{u1.pdf} }
		\thicklines
		\color{darkgreen}
		\put(-17,51){\vector(-1,-1){12.3}}		
		\put(-66,14){\vector(1,0){52}}
		\put(-40,34){$c_1$}
		\put(-75.5,9){$c_2$}
		\color{black}
		\put(-63,32){{$\underline{1}$}}
		\put(-42,5){{$\underline{2}$}}
\end{tabular}
\ \ \ \ \ \ \ \ 
\begin{tabular}{c}
\ \\
		\includegraphics[scale=1.5]{u0.pdf}
		\thicklines
		\color{darkgreen}
		\put(-40.2,17){\vector(1,0){34.5}}
		\put(-50.8,13){$c_2$}
		\color{red}
		\put(-51.2,35){$c_1$}
		\put(-5.7,37){\vector(-1,0){34.5}}
		\color{black}
		\put(-26,43){{$\underline{1}$}}
		\put(-26,24){{$\underline{2}$}}
		\put(-26,6){{$\underline{3}$}}
\ \\
\ \\
\ \\
\ \\
		\includegraphics[scale=1.5]{u00.pdf}
		\thicklines
		\color{darkgreen}
		\put(-26,39){\vector(1,1){12.3}}
		\put(-36.5,34){$c_1$}
		\color{red}
		\put(-78,34){$c_2$}
		\put(-92,51.5){\vector(1,-1){12.3}}
		\color{black}
		\put(-56,11){{$\underline{1}$}}

\end{tabular}
\put(-152,16){\vector(3,1){65}}
\put(-153,-16){\vector(2,-1){38}}
\put(-2,-37){\vector(2,1){51}}
\put(-29.5,38){\vector(3,-1){78}}
\put(-33,36.5){$\circ$}
\put(-224,82){$\overbrace{\ \ \ \ \ \ \ \ \ \ \ \ \ \ \ \ \ \ \ }$}
\put(-91,82){$\overbrace{\ \ \ \ \ \ \ \ \ \ \ \ \ \ \ \ \ \ \ }$}
\put(58,82){$\overbrace{\ \ \ \ \ \ \ \ \ \ \ \ \ \ \ \ \ \ \ }$}
\put(-206,92){$I = -1$}
\put(-72,92){$I = 0$}
\put(77.5,92){$I = 1$}
\ \ \ \ \ \ \ \ \ \ \ \ \ \ \ 
\begin{tabular}{c}
\ \\
		\includegraphics[scale=1.5]{u1.pdf}
		\thicklines
		\color{red}
		\put(-26.2,39){\vector(1,1){12.3}}
		\put(-11,14){\vector(-1,0){52}}
		\put(-37,34){$c_2$}
		\put(-71.5,9){$c_1$}
		\color{black}
		\put(-60,32){{$\underline{1}$}}
		\put(-39,5){{$\underline{2}$}}
\end{tabular}
\end{center}

\begin{center}
\ \ \ \ \ \ \begin{tabular}{c}
\ \\
		\reflectbox{ \includegraphics[scale=1.5]{u1.pdf} }
		\thicklines
		\color{darkgreen}
		\put(-17,51){\vector(-1,-1){12.3}}		
		\put(-14,14){\vector(-1,0){52}}
		\put(-40,34){$c_1$}
		\put(-75,9){$c_2$}
		\color{black}
		\put(-63,32){{$\underline{1}$}}
		\put(-42,5){{$\underline{2}$}}
\end{tabular}
\ \ \ \ \ \ \ \ 
\begin{tabular}{c}
\ \\
		\includegraphics[scale=1.5]{u0.pdf}
		\thicklines
		\color{darkgreen}
		\put(-40.2,17){\vector(1,0){34.5}}
		\put(-50.2,13){$c_2$}
		\color{red}
		\put(-51.2,35){$c_1$}
		\put(-40.2,37){\vector(1,0){34.5}}
		\color{black}
		\put(-26,43){{$\underline{1}$}}
		\put(-26,24){{$\underline{2}$}}
		\put(-26,6){{$\underline{3}$}}
\ \\
\ \\
\ \\
\ \\
		\includegraphics[scale=1.5]{u00.pdf}
		\thicklines
		\color{darkgreen}
		\put(-26,39){\vector(1,1){12.3}}
		\put(-36,34){$c_1$}
		\color{red}
		\put(-78,34){$c_2$}
		\put(-80,39.5){\vector(-1,1){12.3}}
		\color{black}
		\put(-56,11){{$\underline{1}$}}

\end{tabular}
\put(-152,16){\vector(3,1){65}}
\put(-153,-16){\vector(2,-1){38}}
\put(-2,-37){\vector(2,1){51}}
\put(-29.5,38){\vector(3,-1){78}}
\put(-33,36){$\circ$}
\put(-233,-132)
{
	\begin{tabular}{c}
		$\underbrace{\ \ \ \ \ \ \ \ \ \ \ \ \ \ \ \ \ \ \ }$\\
		ker/im\\
		$\shortparallel$\\
		$\vdots$\\
		$\shortparallel$\\
		$\langle \emptyset \rangle$\\
		$\oplus$\\
	 	$\langle \emptyset \rangle$
	\end{tabular}
}
\put(-145,-132)
{
	\begin{tabular}{c}
		$\underbrace{\ \ \ \ \ \ \ \ \ \ \ \ \ \ \ \ \ \ \ }$\\
		ker/im\\
		$\shortparallel$\\
		$\vdots$\\
		$\shortparallel$\\
		$\langle \1 \otimes \1 \otimes x - \1 \otimes x \otimes \1, \  \1 \otimes x \otimes x \rangle$\\
		$\oplus$\\
		$\langle \1 \otimes \1 \otimes x - \1 \otimes x \otimes \1,  \  \1 \otimes x \otimes x \rangle$
	\end{tabular}
}
\put(42,-132)
{
	\begin{tabular}{c}
		$\ \ \underbrace{\ \ \ \ \ \ \ \ \ \ \ \ \ \ \ \ \ \ \ }$\\
		\ \ ker/im\\
		$\ \ \shortparallel$\\		
		\ \ $\vdots$\\
		$\ \ \shortparallel$\\
		\ \ $\langle \emptyset \rangle$\\
		\ \ $\oplus$\\
	 	\ \ $\langle \emptyset \rangle$
	\end{tabular}
}
\ \ \ \ \ \ \ \ \ \ \ \ \ \ \ 
\begin{tabular}{c}
\ \\
		\includegraphics[scale=1.5]{u1.pdf}
		\thicklines
		\color{red}
		\put(-26.2,39){\vector(1,1){12.3}}
		\put(-63,14){\vector(1,0){52}}
		\put(-36.5,34){$c_2$}
		\put(-70.5,9){$c_1$}
		\color{black}
		\put(-60,32){{$\underline{1}$}}
		\put(-39,5){{$\underline{2}$}}
\end{tabular}
\end{center}

\newpage
The $B = 2$ subcomplex:

\begin{center}
\ \ \ \ \ \ \ \ \ \ \ \  \begin{tabular}{c}
\ \\
		\includegraphics[scale=1.5]{u0.pdf}
		\thicklines
		\color{darkgreen}
		\put(-40.2,17){\vector(1,0){34.5}}
		\put(-40.2,37){\vector(1,0){34.5}}
		\put(-51,35){$c_1$}
		\put(-50.6,13){$c_2$}
		\color{black}
		\put(-26,43){{$\underline{1}$}}
		\put(-26,24){{$\underline{2}$}}
		\put(-26,6){{$\underline{3}$}}
\end{tabular}
\ \ \ \ \ \ \ \ \ \ \ \ \ \ \ 
\begin{tabular}{c}
\ \\
		\includegraphics[scale=1.5]{u1.pdf} 
		\thicklines
		\color{red}
		\put(-26.2,39){\vector(1,1){12.3}}
		\put(-35.7,34.5){$c_1$}
		\color{darkgreen}
		\put(-71,9){$c_2$}
		\put(-63,14){\vector(1,0){52}}
		\color{black}
		\put(-60,32){{$\underline{1}$}}
		\put(-39,5){{$\underline{2}$}}
\ \\
\ \\
\ \\
\ \\
		\includegraphics[scale=1.5]{u1.pdf} 
		\thicklines
		\color{red}
		\put(-26.2,39){\vector(1,1){12.3}}
		\put(-36.5,34.5){$c_2$}
		\color{darkgreen}
		\put(-71,9.3){$c_1$}
		\put(-63,14){\vector(1,0){52}}
		\color{black}
		\put(-60,32){{$\underline{1}$}}
		\put(-39,5){{$\underline{2}$}}
\end{tabular}
\put(-138,15){\vector(3,1){54}}
\put(-138,-21){\vector(3,-1){53}}
\put(-1.5,-37.5){\vector(2,1){40}}
\put(-2,37){\vector(2,-1){39}}
\put(-5,35){$\circ$}
\put(-198,82){$\overbrace{\ \ \ \ \ \ \ \ \ \ \ \ \ \ \ \ \ \ \ }$}
\put(-74,82){$\overbrace{\ \ \ \ \ \ \ \ \ \ \ \ \ \ \ \ \ \ \ }$}
\put(60,82){$\overbrace{\ \ \ \ \ \ \ \ \ \ \ \ \ \ \ \ \ \ \ }$}
\put(-179,92){$I = 0$}
\put(-55,92){$I = 1$}
\put(80,92){$I = 2$}
\ \ \ \ \ \ \ \ \ \ \ 
\begin{tabular}{c}
\ \\
		\includegraphics[scale=1.5]{u2.pdf} \ 
		\thicklines
		\color{red}
		\put(-29.5,39){\vector(1,1){12.3}}
		\put(-83.5,39){\vector(-1,1){12.3}}
		\put(-40,34){$c_2$}
		\put(-81,34){$c_1$}
		\color{black}
		\put(-59,10){{$\underline{1}$}}
\end{tabular}
\end{center}

\begin{center}
\ \ \ \ \ \ \ \ \ \ \ \ \begin{tabular}{c}
\ \\
		\includegraphics[scale=1.5]{u0.pdf}
		\thicklines
		\color{darkgreen}
		\put(-40.2,37){\vector(1,0){34.5}}
		\put(-5.5,17){\vector(-1,0){34.5}}
		\put(-51,35.3){$c_1$}
		\put(-50.2,13.3){$c_2$}
		\color{black}
		\put(-26,43){{$\underline{1}$}}
		\put(-26,24){{$\underline{2}$}}
		\put(-26,6){{$\underline{3}$}}
\end{tabular}
\ \ \ \ \ \ \ \ \ \ \ \ \ \ \ 
\begin{tabular}{c}
\ \\
		\includegraphics[scale=1.5]{u1.pdf} 
		\thicklines
		\color{red}
		\put(-26.2,39){\vector(1,1){12.3}}
		\put(-37,34.5){$c_1$}
		\color{darkgreen}
		\put(-71.6,9){$c_2$}
		\put(-10.5,14){\vector(-1,0){52}}
		\color{black}
		\put(-60,32){{$\underline{1}$}}
		\put(-39,5){{$\underline{2}$}}
\ \\
\ \\
\ \\
\ \\
		\includegraphics[scale=1.5]{u1.pdf} 
		\thicklines
		\color{red}
		\put(-26.2,39){\vector(1,1){12.3}}
		\put(-36.5,34.5){$c_2$}
		\color{darkgreen}
		\put(-71,9){$c_1$}
		\put(-10.5,14){\vector(-1,0){52}}
		\color{black}
		\put(-60,32){{$\underline{1}$}}
		\put(-39,5){{$\underline{2}$}}
\end{tabular}
\put(-138,15){\vector(3,1){54}}
\put(-138,-21){\vector(3,-1){53}}
\put(-1.5,-37.5){\vector(2,1){40}}
\put(-2,37){\vector(2,-1){39}}
\put(-5,35){$\circ$}
\put(-272,-132)
{
	\begin{tabular}{c}
		$\underbrace{\ \ \ \ \ \ \ \ \ \ \ \ \ \ \ \ \ \ \ }$ \ \ \ \ \ \\
		ker/im \ \ \ \ \ \\
		$\shortparallel$ \ \ \ \ \ \\
		$\vdots$ \ \ \ \ \ \\
		$\shortparallel$ \ \ \ \ \ \\
		$\langle \1 \otimes x \otimes x -  x \otimes \1 \otimes x +  x \otimes x \otimes \1, \ x \otimes x \otimes x \rangle \ $\\
		$\oplus$ \ \ \ \ \ \\
		$\langle \1 \otimes x \otimes x -  x \otimes \1 \otimes x +  x \otimes x \otimes \1, \ x \otimes x \otimes x \rangle \ $
	\end{tabular}
}
\put(-83.5,-132)
{
	\begin{tabular}{c}
		$\underbrace{\ \ \ \ \ \ \ \ \ \ \ \ \ \ \ \ \ \ \ }$\\
		\ ker/im\\
		\ \ $\shortparallel$\\
		\ \ $\vdots$\\
		\ \ $\shortparallel$\\
		\ \ $\langle \emptyset \rangle$\\
		\ \ $\oplus$\\
	 	\ \ $\langle \emptyset \rangle$
	\end{tabular}
}
\put(44,-132)
{
	\begin{tabular}{c}
		$\ \ \underbrace{\ \ \ \ \ \ \ \ \ \ \ \ \ \ \ \ \ \ \ }$\\
		\ \ ker/im\\
		$\ \ \shortparallel$\\		
		\ \ $\vdots$\\
		$\ \ \shortparallel$\\
		\ \ $\langle \emptyset \rangle$\\
		\ \ $\oplus$\\
	 	\ \ $\langle \emptyset \rangle$
	\end{tabular}
}
\ \ \ \ \ \ \ \ \ \ \ 
\begin{tabular}{c}
\ \\
		\includegraphics[scale=1.5]{u2.pdf} \ 
		\thicklines
		\color{red}
		\put(-29.5,39){\vector(1,1){12.3}}
		\put(-39.7,34.5){$c_2$}
		\put(-80.7,34){$c_1$}
		\put(-95.5,51.2){\vector(1,-1){12.3}}
		\color{black}
		\put(-59,10){{$\underline{1}$}}
\end{tabular}
\end{center}

The above calculations of the diagramless homology for the unknot are summarized below.  The copies of $\Z$ are grouped to show the correspondence to Khovanov homology -- two copies of $(\Z \oplus \Z)$ in each $B$-grading.

$$
\D^{i,b}_2 ([ \bigcirc ]) \cong  
	\left\{ \begin{array}{cc} 
		(\Z \oplus \Z) \oplus (\Z \oplus \Z) & \textrm{for } \ i=0, \ b=0\\ 
		(\Z \oplus \Z) \oplus (\Z \oplus \Z) & \textrm{for } \ i=0, \ b=1\\ 
		(\Z \oplus \Z) \oplus (\Z \oplus \Z) & \textrm{for } \ i=0, \ b=2\\ 
		0 & \textrm{otherwise} 
	\end{array} \right .
$$

\section{Constructing Diagramless Homology in Other 3-manifolds}

One of the main benefits of realizing the Khovanov homology of a link using this diagramless approach is its potential to be generalized to links in 3-manifolds other than the 3-sphere.  In Section \ref{2} we viewed the surface $F$ as being contained in $S^3$.  Could we have used 3-manifolds other than $S^3$?  The answer is yes.  

\subsection{The Homological Grading}

One issue is that the homological grading is defined in terms of linking number.   It is well known that the definition of linking number can be extended to simple closed curves in a homology sphere.  In this way the same homological grading could be used for $D^k$-surfaces in a homology sphere.  Defining a version of linking number in more general manifolds has been studied recently by Chernov and Rudyak (see 
\cite{C-R}).  

Another way to define the homological grading for surfaces in 3-manifolds would be to use a local grading.  That is, simply let $I = (\#$ of \textcolor{darkgreen}{active} crosscuts of $F)$; this results in a bonafide homological grading.  The difference is that the $N$ copies of Khovanov homology that appeared before all had the same signature, but the underlying surfaces had different numbers of \textcolor{darkgreen}{active} crosscuts.  Hence, corresponding copies of homology will no longer have the same homological grading.

\subsection{The Definition of $D^k$-surface}

Another potential issue is how one should define $D^k$-surface in a 3-manifold other than $S^3$.  Recall that $S^3 = B^3_+ \cup B^3_-$, and that the definition of $D^k$-surface requires the surface to be able to be isotoped so that the skeleton of the cross-dual is embedded in $\Sigma = B^3_+ \cap B^3_-$ and the facets are properly embedded in $B^3_+$. 

Translating this condition from $S^3$ to another 3-manifold $M$ requires $M$ to have a (fixed) Heegaard splitting.  Let $A$ and $B$ be handlebodies such that $M = A \cup B$, with Heegaard surface $S = A \cap B$.  Then the following definition of $D^k$-surface is the proposed one for links in arbitrary closed oriented 3-manifolds.

\begin{defn}
A {\it $D^k$-surface} is a compact surface $F \subseteq M = A \cup B$ with $k$ crosscuts $\{c_1, ...,c_k\}$ such that
\begin{itemize}
\item the crosscuts are oriented and ordered,
\item the facets of $F$ are allowed to be decorated by dots (which are \textit{not} allowed to move from one facet to another), 
\item the cross-dual $F^{\textrm{cd}}$ is orientable (this global orientatibility is independent of the local orientations of the cross-dual),
\item skel$(\fcd) \subseteq S$ with all of the locally oriented pieces of surface of $\fcd$ agreeing with the orientation of $S$, and
\item $\fcd - \textrm{skel}(\fcd) = \{\textrm{the facets of }\fcd\} \subseteq A$.
\end{itemize}
\end{defn}

\subsection{Computing Examples}

Although this diagramless homology theory might be easily defined for links in any closed orientable 3-manifold $M$, computing examples could be difficult.  In this paper, the $D^k$-surfaces involved in an example are easy to consider because of the ease of visualizing surfaces in 3-space.  In addition to this, Proposition \ref{big_prop} allows us to work with state surfaces instead of $D^k$-surfaces; it is unlikely that such a proposition can be exploited for links in more general 3-manifolds.  Other methods will have to be used to successfully compute examples.

\section{Additional Remarks}
\label{70}

\subsection{Using Boundary Slope as the Homological Grading}

A possible improvement to this theory would be the ability to distinguish between different orientations of links.  A link's orientation does not factor in when calculating the signature of a surface.  However, there is an alternative to surface signature which does detect link orientation: boundary slope.  

The \textit{boundary slope} of a surface is the linking number of the boundary of the surface with the pushoff of the boundary in the tangential direction away from the surface.  Although signature seems to be the natural choice for the homological grading for the diagramless homology theory, the following proposition would allow the use of boundary slope as the homological grading, hence distinguishing between different orientations of a given link.

\begin{prop}
The differential $d$ for the diagramless homology theory increases boundary slope by +2.  Therefore, for a collection surfaces in the same diagramless subcomplex, the boundary slope and two times the signature differ by a constant.
\end{prop}

\emph{Proof.}  It is straightforward to calculate that $d_c$, hence $d$, increases boundary slope by +2.  In the figure given below, one can calculate that the linking number of the blue and red pushoffs with the boundary are (locally) 0 on the left and (locally) +2 on the right.  Reversing the orientation of one or both of the red or blue lines (as well as the corresponding parts of the boundary) gives the same result.

\begin{center}
\includegraphics[scale=2]{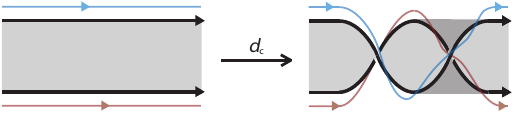}
\end{center}  \qed

\subsection{The Value of \texorpdfstring{$N$}{\textit{N}} in Theorem \ref{decomp_thm}}

Theorem \ref{decomp_thm} states that, for a fixed number of crosscuts, $k$, the diagramless homology of a link is equal to the direct sum of $N$ copies of the Khovanov homology of that link.  The value of $N$ could determine a link invariant if, for example, we let $k$ equal the minimum crossing number of the link.  The particular value of $N$ for a given number of crosscuts $k$ is not explored in depth in this paper.  However, two conjectures are given concerning the value of $N$.

\begin{conj}
If the number of crosscuts, $k$, is less than the minimum crossing number of a link, then $N = 0$.
\end{conj}

The motivation for the next conjecture comes from Definition \ref{e_class}, where the link diagram equivalence classes are defined.  Essentially, two diagrams (for a non-split link are equivalent if they are isotopic in 2-space, after possibly turning one diagrams upside-down (flipping it).  This extra flipping means that the usual number of link diagram equivalence classes (under 2-space isotopy only) could differ from the number of our equivalence classes by up to a factor of two.  

\begin{conj}
Let $L$ be a non-split link and $k$ be a fixed positive integer.  If $n_k(L)$ is the number of distinct $k$-crossing link diagrams for $L$ up to 2-space isotopy, then 
$$
\frac{n_k(L)}{2} \leq N \leq n_k(L),
$$
where $N$ is the number from Theorem \ref{decomp_thm}, which equals the number of copies of Khovanov homology in the diagramless homology of $L$ with $k$ crosscuts.
\end{conj}

If the above conjecture is true, we then have the following corollary.  However, this corollary could likely be proved without the above conjecture.

\begin{cor}
Let $L$ be a link.  If there is only one $k$-crossing diagram for $L$ up to 2-space isotopy, then $N = 1$.
\end{cor}

Putting tight bounds on $N$ when $L$ is a split link may be more difficult.  The way in which link diagram equivalence classes are defined in this paper allows for disconnected components of diagrams to be moved around and be embedded in different regions of the other diagram component(s).  

\subsection{Morphisms Between Diagramless Complexes}

In one of his papers on Khovanov homology (\cite{BN2}), Bar-Natan uses smoothings of enhanced Kauffman states of link diagrams as objects in a category he calls $\mathcal{C}ob^3$.  In this category, morphisms are the cobordisms between such smoothings.  Below we explore possibility to define a similar category for the diagramless theory.

At first, it would seem that the corresponding category for the diagramless theory would have $D^k$-surfaces as the objects instead of diagram smoothings.  However, Proposition \ref{big_prop} allows us to work with state surfaces instead of $D^k$-surfaces, and state surfaces are built from smoothings of enhanced Kauffman states.  Hence, the objects in the category for the diagramless theory might be represented by smoothings as well.  In this case the morphisms would be represented by cobordisms.  

In \cite{C-S}, Carter and Saito use movie moves to study such cobordisms.  Carter and Saito introduce movie moves by showing them alongside their corresponding cobordisms.  Since we would be considering the state surfaces corresponding to the smoothings involved, we would be interested in `movies of state surfaces'.  A few examples are given below.
\bigskip
\begin{center}
\includegraphics[trim = 0mm 107mm 0mm 0mm, clip,scale=1]{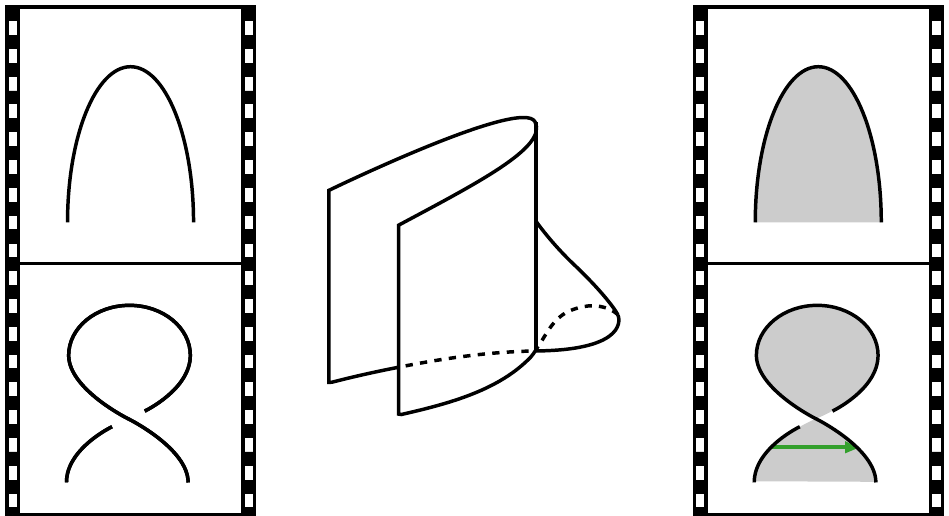}\\
{\footnotesize The Reidemeister 1 move}\\
\ 
\end{center}

\begin{center}
\includegraphics[trim = 0mm 53.5mm 0mm 53.5mm, clip,scale=1]{movies.pdf}\\
{\footnotesize Deletion of a circle // a death (viewed from bottom to top)}\\
\ 
\end{center}

\begin{center}
\includegraphics[trim = 0mm 0mm 0mm 107mm, clip,scale=1]{movies.pdf}\\
{\footnotesize A smoothing change near a crossing // a saddle}\\ \ 
\end{center}

Due to the presence of crosscuts which are ordered, oriented, and given a label of \textcolor{green}{active} or \textcolor{red}{inactive}, additional information must be given along with the movie or cobordism representing the morphism at hand.  Determining how best to do this requires more work, but a thorough treatment of morphisms is not given in this paper.

\end{document}